\documentclass[11pt]{article}
\usepackage[margin=0.75in]{geometry}

\usepackage{amssymb,amsmath,mathtools,amsthm,thmtools,amsfonts,dsfont,mathrsfs}
\usepackage{thm-restate}
\usepackage[normalem]{ulem}

\usepackage{url} 
\usepackage{hyperref}
\hypersetup{
colorlinks,
linkcolor={red}, 
citecolor={black},
urlcolor={blue!60!black},
pdftitle={Hypergraph Hajnal-Szemeredi},
pdfauthor={Michelle Delcourt, Cece Henderson, Thomas Lesgourgues, Luke Postle}}

    \usepackage{listings}
    \usepackage{xcolor}

    \definecolor{codegreen}{rgb}{0,0.6,0}
    \definecolor{codegray}{rgb}{0.5,0.5,0.5}
    \definecolor{codepurple}{rgb}{0.58,0,0.82}
    \definecolor{backcolour}{rgb}{0.95,0.95,0.92}

    \lstdefinestyle{mystyle}{
        backgroundcolor=\color{backcolour},
        commentstyle=\color{codegreen},
        keywordstyle=\color{magenta},
        numberstyle=\tiny\color{codegray},
        stringstyle=\color{codepurple},
        basicstyle=\ttfamily\footnotesize,
        breakatwhitespace=false,
        breaklines=true,
        captionpos=b,
        keepspaces=true,
        numbers=left,
        numbersep=5pt,
        showspaces=false,
        showstringspaces=false,
        showtabs=false,
        tabsize=2
    }
    
\usepackage[utf8]{inputenc}
\usepackage{color}
\usepackage{float} %for figure placement
\usepackage{subcaption} %for subfigure
\usepackage{enumitem}
\setlist[itemize]{topsep=.1in,itemsep=.1in}
\setlist[enumerate]{topsep=.1in,itemsep=.1in}
\setlist[enumerate,1]{label={(\roman*)}}
\setlist[enumerate,2]{label={(\alph*)}}
\setlist[enumerate,3]{label={(\arabic*)}}

\setlist[itemize]{nolistsep,noitemsep, topsep=0pt}

\newcommand{\eps}{\varepsilon}

\newcommand{\mc}[1]{\mathcal{#1}}%

\DeclarePairedDelimiter{\parens}{(}{)}
\DeclarePairedDelimiter{\set}{\{}{\}}
\DeclarePairedDelimiter{\brackets}{[}{]}
\DeclarePairedDelimiter{\floor}{\lfloor}{\rfloor}
\DeclarePairedDelimiter{\ceil}{\lceil}{\rceil}
\DeclarePairedDelimiter{\abs}{\lvert}{\rvert}   % |absolute value|
\DeclarePairedDelimiter{\size}{\lvert}{\rvert}   % |size|

\usepackage[noabbrev,capitalise]{cleveref}

\theoremstyle{plain}
\newtheorem{thm}{Theorem}[section]
\newtheorem{lem}[thm]{Lemma}

\newtheorem{cor}[thm]{Corollary}

\newtheorem{question}[thm]{Question}

\newenvironment{proofclaim}[1][\proofname]
{\proof[#1]}
{\endproof}

\declaretheorem[
  style=plain,
  name=Claim,
  within=theorem,
]{claim}

\newenvironment{lateproof}[1]
 {%
  \renewcommand{\theclaim}{\ref{#1}.\arabic{claim}}%
  \begin{proof}[Proof of~\cref{#1}]%
 }
 {\end{proof}}

\usepackage{etoolbox}
 \AtEndEnvironment{proof}{\setcounter{claim}{0}}

\theoremstyle{plain} % just in case the style had changed
\newcommand{\thistheoremname}{}
\newtheorem{genericthm}{\thistheoremname}

\theoremstyle{definition}
\newtheorem{definition}[thm]{Definition}

\crefname{equation}{}{}
\crefname{lem}{Lemma}{Lemmas}
\crefname{claim}{Claim}{Claims}
\crefname{thm}{Theorem}{Theorems}
\crefname{enumi}{}{}
\crefname{lstlisting}{Code}{Codes}

\usepackage[dvipsnames, table]{xcolor}

\title{Hypergraph Hajnal-Szemer\'edi}

\author{
Michelle Delcourt
\thanks{Department of Mathematics, Toronto Metropolitan University (formerly named Ryerson University),
Toronto, Ontario M5B 2K3, Canada \texttt{mdelcourt@torontomu.ca}. Research supported by NSERC under Discovery Grant No. 2019-04269 and a Sloan Research Fellowship.} 
\and 
Cicely (Cece) Henderson
\thanks{Combinatorics and Optimization Department,
University of Waterloo, Waterloo, Ontario N2L 3G1, Canada 
\{\texttt{c3hender,lpostle}\}\texttt{@uwaterloo.ca}. Partially supported by NSERC
under Discovery Grant No. 2026-04610.} \and
Thomas Lesgourgues \thanks{School of Mathematics and Statistics, UNSW Sydney, NSW 2052, Australia \texttt{t.lesgourgues@unsw.edu.au}. Partially supported by Australian Research Council grant DP250101611.}
\and
Luke Postle \footnotemark[2] 
}
\date{\today}

\begin{document}

\maketitle

\begin{abstract}
A seminal result of Hajnal and Szemer\'edi from 1970 asserts that every graph on $n$ vertices with $q\;|\;n$ and minimum degree at least $\Big(1-\frac{1}{q}\Big)n$ admits a $K_q$-factor (a collection of vertex-disjoint copies of $K_q$ spanning all vertices). There has been sustained interest in extending this result to hypergraph $K_q^r$-factors. For general $q$ and $r$, Lu and Sz\'ekely in 2007 showed that every large enough, $r$-graph $G$ with $q\mid v(G)$ and $\delta_{r-1}(G) \ge \left(1- \frac{c}{q^{r-1}}\right)v(G)$ admits a $K_q^r$-factor (where $c$ is a constant depending on $r$ but not $q$), and there have been various subsequent improvements to the constant (but not the order) of this bound over the ensuing decades. Meanwhile, results from codegree Tur\'an theory show that $\delta_{r-1}(G) \ge \left(1 - \frac{c \cdot \log q}{q^{r-1}} \right)v(G)$ is necessary.

Surprisingly, our main result shows that the order of this upper bound is not correct and, in fact, the codegree $K_q^r$-factor threshold is much closer to this codegree Tur\'an theory bound. 
In particular, we show that $\delta_{r-1}(G) \ge \left( 1-c\cdot \frac{\log^{1-1/(r-1)}q}{q^{r-1}}\right)v(G)$ suffices for the existence of a $K_q^r$-factor. 
Our key contribution is showing the same is true for the fractional threshold. Inspired by the first and fourth author's recent proof of the Nash-Williams' Conjecture, we accomplish this by combining a linear programming framework with ideas from discharging.
Additionally, we establish a technical result about finding large independent sets in sparse, non-uniform hypergraphs, which may be of independent interest. 
\end{abstract}

%%%%%%%%%%%%%%%%%%%%%%
% Introduction
%%%%%%%%%%%%%%%%%%%%%%

\section{Introduction}

A central thread in extremal combinatorics is the question of density thresholds for the existence of certain structures in graphs. A primary example of such structures is the triangle, and a classical result, often attributed to Mantel~\cite{mantel} in 1907, asserts that every graph on $n$ vertices with at least $\lfloor \frac{n^2}{4} \rfloor +1$ edges contains $K_3$ as a subgraph. In 1941, Tur{\'a}n~\cite{turan1941extremal} generalized this result to larger cliques by determining the exact number of edges needed to guarantee the existence of a copy of $K_q$ in a graph - roughly around $\big(1-\frac{1}{q-1}\big)n$ - and complete balanced multipartite graphs show this is tight. This became known as Tur\'an's theorem, the numbers as Tur\'an numbers, and launched an entire subfield at the heart of extremal combinatorics known as Tur\'an theory.

Tur\'an's question about the existence of a $q$-clique in a graph $G$ is equivalent to asking whether the empty set of vertices is contained in a copy of $K_q$ in $G$. Thus, it is natural to generalize this to the question of the existence of a collection of cliques containing every size 1 subset of vertices (i.e. every vertex). In particular, a \emph{$K_q$-factor} is a set of vertex-disjoint copies of $K_q$ spanning all vertices. Here, minimum degree is a natural measure of density since an isolated vertex is not contained in any copy of $K_q$. Also note that $q \mid v(G)$ is a necessary divisibility condition for the existence of a $K_q$-factor in $G$. For the case of triangles, the famous Corr{\'a}di-Hajnal~\cite{CH63} Theorem of 1963 shows that any graph $G$ on $n$ vertices with $3 ~|~n$ and $\delta(G) \ge \frac{2}{3} n$ contains a $K_3$-factor. The following seminal result of Hajnal and Szemer\'edi~\cite{HS70} from 1970 generalizes this to $K_q$-factors:

\begin{thm}[Hajnal-Szemer\'edi \cite{HS70}]\label{thm:HS}
	For every integer $q\geq 3$, every graph $G$ on $n$ vertices with $q \mid n$ and
    \[\delta(G)\geq \Big(1-\frac{1}{q}\Big)n\]
    contains a $K_q$-factor.
\end{thm}
We note that Kierstead and Kostochka~\cite{KK08} provided a proof of~\cref{thm:HS} in 2008 that is not only significantly shorter than the original but also more directly phrased in terms of clique factors. Hajnal and Szemer\'edi's argument is written for equitable colorings, as we discuss in~\cref{Subsec:Connection}.
This minimum degree threshold is best possible: consider a complete $q$-partite graph $G$ with one part of size $\ceil{n/q}+1$ and others as equal as possible. In this case, $\delta(G)=(1-\frac{1}{q})n-1$, yet $G$ contains no $K_q$-factor.

The purpose of this paper is to make substantial progress on a hypergraph Hajnal-Szemer\'edi theorem.

\subsection{Main Results}\label{Subsec:IntroFactors}
Given the centrality and extensive history of the above questions, it is natural to consider extending them to the case of complete hypergraphs. Let $K_q^r$ denote the complete $r$-uniform hypergraph on $q$ vertices. For brevity, we say an $r$-uniform graph $G$ is an \emph{$r$-graph}. 

In terms of the existence of one copy of $K_q^r,$ the problem of determining hypergraph Tur\'an numbers is notoriously difficult. For general $q$ and $r$, a 1981 lower bound of Sidorenko~\cite{S81} and a 1983 upper bound of de Caen~\cite{dC83} together imply that the Tur\'an number for $K_q^r$ is $1 - \Theta_r\left(\frac{1}{q^{r-1}} \right)$, thus settling the order for fixed $r$ and large $q$. For general $q,$ in 2025 Pikhurko~\cite{P25} disproved a longstanding conjecture of de Caen by showing that, for $q = r+s$, the trivial lower bound of $1 - \frac{1}{\binom{r+s}{s}}$ is tight up to a multiplicative constant depending only on $s$.

What about a hypergraph version of the Hajnal-Szemer\'edi theorem? Such a result would not only be a natural extension of~\cref{thm:HS}, a classical theorem at the core of extremal combinatorics, but also it would have connections to other important concepts throughout the field. In particular, we detail the connection between a hypergraph Hajnal-Szemer\'edi theorem to equitable coloring, powers of Hamilton cycles, and combinatorial designs in~\cref{Subsec:Connection}. Additionally, there has been sustained interest in generalizing this fundamental result to other settings, such as Keevash and Mycroft's~\cite{KM15} 2015 multipartite Hajnal-Szemer\'edi theorem and Treglown's~\cite{T15} 2015 directed Hajnal-Szemer\'edi theorem. However, a hypergraph version has remained elusive.  
As with many graph theoretic concepts, there are multiple ways to extend the notion of minimum degree to hypergraphs. One of the most natural generalizations that has garnered the most attention is known as \emph{minimum codegree} or \emph{minimum $(r-1)$-degree}: let $\delta_{r-1}(G)$ denote the minimum over all $(r-1)$-sets $S \subseteq V(G)$ of the number of edges of $G$ containing $S$. Codegree has been employed for notable results in the related area of design theory, such as in all proofs of the infamous Existence Conjecture for Combinatorial Designs~\cite{RefineI, GKLO16,  K24, K14}.
Thus, we are motivated to consider the following question.

\begin{question}[Hypergraph Hajnal-Szemer\'edi Question]\label{HHSQuestion}
What is the minimum codegree threshold for guaranteeing the existence of a $K_q^r$-factor?
\end{question}

As we will discuss, the order of this threshold for fixed $r$ and large $q$ is still open. The main result of this paper is a significantly improved upper bound, thus making a step towards determining the correct order of this minimum codegree threshold.

Exact thresholds for hypergraph clique factors are known for a few cases of $q$ and $r$. For example, in 2009 R\"{o}dl, Ruci\'{n}ski, and Szemer\'edi~\cite{RRS09} found the threshold for the existence of a perfect matching in an $r$-graph, Keevash and Mycroft~\cite{KM15} determined the threshold for a $K_4^3$-factor in 2015, and Falgas-Ravry and Zhao~\cite{FZ16} in 2016 as well as Han, Lo, Treglown, and Zhao~\cite{HLTZ17} in 2017 investigated minimum codegree thresholds for other small 3-uniform graphs. We refer readers to the 2016 survey of Zhao~\cite{Z16} for a more comprehensive history.

For general $q$ and $r$, the best known lower bound for~\cref{HHSQuestion} comes from codegree Tur\'an theory. 
The \textit{codegree Tur\'an density} for $K_q^r$, denoted $T_{\text{codeg}}(K_q^r)$, is the minimum $\gamma >0$ such that any $r$-graph on $n$ vertices with $\delta_{r-1}(G) \ge (\gamma + o(1))\cdot n$ contains $K_q^r$ as a subgraph. As noted by Glock, K{\"u}hn, Lo, and Osthus in \cite{GKLO16}, it follows from a 2014 construction of Kostochka, Mubayi, and Verstra{\"e}te~\cite{KMV14} that there exist infinitely many $K_q^r$-free $r$-graphs $G$ with $\delta_{r-1}(G) \geq \left(1 - c\cdot \frac{\log q}{q^{r-1}}\right)\cdot v(G)$ for a constant $c$ depending only on $r$. 
In 2018 Lo and Zhao~\cite{LZ2018codegree} derived an upper bound of the same order, thus determining $T_{\text{codeg}}(K_q^r)$ as follows. 
\begin{thm}[Lo and Zhao~\cite{LZ2018codegree}]\label{Thm:LB}\label{thm:HyperTuran}
For every integers $3 \le r< q$,
    \[T_{{\rm codeg}}(K_q^r) = 1- \Theta_r\left(\frac{\log q}{q^{r-1}}\right).\] 
\end{thm}
It is striking in that the $\log q$ factor is needed for the codegree Tur\'an density when it is not for the $K_q^r$ Tur\'an number. In the context of~\cref{HHSQuestion}, \cref{Thm:LB} provides the best known lower bound on the minimum codegree threshold for the existence of a clique factor.  

In 2007, Lu and Sz\'ekely~\cite{LS07} provided the first upper bound for the minimum codegree $K_q^r$-factor threshold for general $q$ and $r$ as follows. 

\begin{thm}[Lu and Sz\'ekely~\cite{LS07}]\label{thm:PreviousHHS}
    For every integer $r\geq 2$, there exists a constant $c > 0$ such that the following holds for every $q>r$. Every large enough $r$-graph $G$ with 
    \[\delta_{r-1}(G) \geq \left(1 - \frac{c}{q^{r-1}}\right)\cdot v(G),\]
    admits a $K_q^r$-factor.
\end{thm}
(We utilize a corollary of this result, \cref{cor:FracHHR}, in our subsequent arguments.) Lu and Sz\'ekely proved~\cref{thm:PreviousHHS} more generally in terms of minimum vertex degree. The minimum codegree bound in~\cref{thm:PreviousHHS} implies the corresponding required minimum vertex degree and hence is a corollary of their result.
The constant on the bound in~\cref{thm:PreviousHHS} was subsequently improved in 2021 by Pavez-Sign{\'e}, Sanhueza-Matamala, and Stein~\cite{PSS21}:
\begin{thm}[Pavez-Sign{\'e}, Sanhueza-Matamala, and Stein~\cite{PSS21}]\label{thm:prevHHSImprovedConst}
    For every integer $r\geq 2$, there exists a constant $c > 0$ such that the following holds for every $q>r$. Every large enough $r$-graph $G$ with
    $$\left( 1 - \frac{1}{\binom{q-1}{r-1}+\binom{q-2}{r-2} }\right)v(G)$$
    admits a $K_q^r$-factor. 
\end{thm}
 We note that Pavez-Sign{\'e}, Sanhueza-Matamala, and Stein's proved~\cref{thm:prevHHSImprovedConst} in the context of a hypergraph version of the Pósa-Seymour conjecture. We discuss the connection between factors and powers of Hamilton cycles in~\cref{Subsec:Connection}. An alternate proof of~\cref{thm:prevHHSImprovedConst} was given by Lang and Sanhueza-Matamala~\cite{LS24} in 2024.

We remark on the large gap between the lower bound for the codegree $K_q^r$-factor threshold given by~\cref{Thm:LB} and the upper bound provided by~\cref{thm:PreviousHHS}. Indeed, in 2015 Lo and Markstr\"{o}m~\cite{LM15} noted this gap and specifically asked whether the $\log  q$ factor is necessary in the three uniform case:
% In 2015, Lo and Markstr\"{o}m~\cite{LM15} developed bounds for 3-uniform cliques 
% %by showing that the packing threshold is between
% %\[\Big(1-\frac{196 \log(q+1)}{(q-1)^2}\Big)n\quad\text{ and }\quad \Big(1-\frac{2}{q^2-3q+4}+o(1)\Big)n.\]
% and, for larger uniformities, they showed that for large enough $r$ and large enough $q$ (depending on $r$), the minimum codegree $K_q^r$-factor threshold is at most $\left(1 - \binom{q-1}{r-1}^{-1} + \eps\right)$. In their paper, Lo and Markstr\"{o}m~\cite{LM15} remark on the gap between this lower and upper bound and specifically asked whether the $\log  q$ factor is necessary in the three uniform case:
\begin{question}[Question 1.7 in~\cite{LM15}]\label{LMQuestion}
    How does the minimum codegree threshold for the existence of a $K_q^r$-factor behave for $q\geq r\geq 3$ as $n\to\infty$? When $r = 3,$ is this threshold equal to $(1-C\log q / q^2)n$ or $(1-C/ q^2)n$ for large $n$ and some absolute constant $C$ independent of $q$ and $n$?
\end{question}

Surprisingly, our main result shows that the order of the upper bound given by~\cref{thm:PreviousHHS} is not correct and, in fact the codegree $K_q^r$-factor threshold is much closer to the codegree Tur\'an theory bound from~\cref{Thm:LB}.

\begin{thm}\label{thm:mainexact}
    For every integer $r \geq 3$, there exists a real $c > 0$ such that for all integers $q > r$, the following holds for every large enough integer $n$ such that $q\mid n$. If $G$ is an $r$-uniform hypergraph on $n$ vertices with
    \[\delta_{r-1}(G)\geq \parens*{1-c\cdot \frac{\log^{1-1/(r-1)}q}{q^{r-1}}}\cdot n,\]
    then $G$ contains a $K_q^r$-factor.
\end{thm}

Therefore, \cref{thm:mainexact} answers~\cref{HHSQuestion} within a factor $\log^{1/(r-1)}q$ and~\cref{LMQuestion} for general $q$, also within a factor of $\log^{1/(r-1)}q$.
We note that this bound nicely matches~\cref{thm:HS} when $r=2$, and thus one might be tempted to conjecture that the minimum codegree $K_q^r$-factor threshold is in fact $1 - \Theta_r\parens*{\frac{\log^{1-1/(r-1)}q}{q^{r-1}}}$. We also note that we have not considered what our techniques give for a $d$-degree threshold more generally. It might be the case that the minimum $d$-degree threshold for a $K_q^r$-factor forms a sequence where the $\log q$ term is raised to the power of $(d-1)/(r-1)$ when considering values of $d$ ranging between $1$ and $r-1$. 

To prove~\cref{thm:mainexact}, we utilize a result of Lang~\cite{lang2023tiling} that ties the minimum codegree threshold for a $K_q^r$-factor to thresholds for avoiding space, divisibility, and covering barriers to the existence of such a factor. In the context of clique factors, the space barrier poses the most significant obstruction and thus~\cref{thm:mainexact} follows from our substantial improvement on the fractional threshold for a $K_q^r$-factor, as detailed in~\cref{Subsec:Barriers}. 

\subsection{Connection to Equitable Coloring, Powers of Hamilton Cycles, and Designs}\label{Subsec:Connection}
Originally, Hajnal and Szemer\'edi phrased~\cref{thm:HS} in terms of equitable colorings. A \emph{$t$-equitable coloring} of a graph $G$ is a proper vertex coloring where every color class has size exactly $\lfloor \frac{v(G)}{t}\rfloor$ or $\lceil \frac{v(G)}{t}\rceil$. By considering graph complements, \cref{thm:HS} is equivalent to saying that a graph $G$ with $\Delta(G) < n/q$ admits an $n/q$-equitable coloring. Equitable coloring is defined analogously for hypergraphs, where a proper coloring is one in which no edge is monochromatic.
Previous generalizations of equitable colorings have focused on bounding the number of edges for the existence of an equitable coloring, such as Akhmejanova and Shabanov's~\cite{AS22} 2022 result, or on providing an equitable coloring vertex degree threshold for small uniformities, such as the 2010 work of Kierstead and Mubayi~\cite{KM10} for $3$-graphs. Alternatively,~\cref{thm:mainexact} provides a codegree threshold for equitable coloring of hypergraphs of all uniformities.

\begin{cor}
    For every integer $r \geq 3$, there exists a real $c > 0$ such that for all integers $q > r$, the following holds for every large enough integer $n$ with $q\mid n$. If $G$ is an $r$-uniform hypergraph on $n$ vertices with
    \[\Delta_{r-1}(G)<c\cdot \frac{\log^{1-1/(r-1)}q}{q^{r-1}}\cdot n,\]
    then $G$ admits an $\frac{n}{q}$-equitable coloring.
\end{cor}
We note that though $q\mid n$ is not required for an equitable coloring, we do still require this divisibility condition to apply~\cref{thm:mainexact}. (This is not the case for the graph setting, where one can add a clique of auxiliary vertices to make a graph of bounded maximum degree divisible, then apply~\cref{thm:HS} to the complement.) Additionally, our proof methods require that $n$ be large enough, which is not necessary in~\cref{thm:HS}.

The clique factors of~\cref{thm:mainexact} are also related to powers of Hamilton cycles. The \emph{$q$th power of a Hamilton cycle} is a Hamilton cycle with an additional chord between every pair of vertices of distance at most $q$ along the cycle. Since the existence of a $K_q^r$-factor is necessary for a $(q-1)$th power of a Hamilton cycle,~\cref{thm:HS} provides a trivial lower bound on the threshold for the existence of powers of Hamilton cycles. In fact, 
%generalizing Dirac's~\cite{dirac1952some} classical result on the minimum degree threshold for a Hamilton cycle, 
P\'osa (for $q = 3$) and Seymour (for $q > 3$) famously conjectured that this minimum degree threshold of $\delta(G) \ge (1 - \frac{1}{q})n$ suffices for the existence of a $(q-1)$th power of a Hamilton cycle, and Koml\'os, S\'ark\"ozy, and Szemer\'edi~\cite{KSS98} confirmed this conjecture in 1998 for all values of $q$ and large enough $n$. In terms of hypergraphs, in 2008 R{\"o}dl, Rucinski, and Szemer{\'e}di~\cite{RRS08} extended Dirac's classical result about minimum degree thresholds for Hamilton cycles to hypergraphs by showing that an $r$-graph on $n$ vertices with minimum codegree at least $(1+o(1))n/2$ contains a tight Hamilton cycle (a cycle with maximum overlap between edges). Recall from~\cref{thm:prevHHSImprovedConst} that in 2021, Pavez-Sign{\'e}, Sanhueza-Matamala, and Stein~\cite{PSS21} showed that minimum codegree at least $\left( 1 - 1/\left(\binom{q-1}{r-1}+\binom{q-2}{r-2} \right)\right)v(G)$ suffices for the existence of a $(q-r+1)$th power of a Hamilton cycle in an $r$-graph $G$. Whether the codegree threshold for a $K_q^r$-factor is the same for that of a $(q-r+1)$th power of a tight Hamilton cycle in $r$-graphs remains an open question, and thus~\cref{thm:main} is relevant to a hypergraph version of the influential P\'osa-Seymour conjecture.

Additionally, \cref{thm:mainexact} is connected to the field of design theory. A \emph{$K_q^r$-decomposition} of an $r$-graph $G$ is a set of edge-disjoint copies of $K_q^r$ spanning the edges of $G$. The celebrated Nash-Williams' conjecture~\cite{nash1970unsolved} of 1970 posits that a graph $G$ on $n$ vertices which satisfies the necessary divisibility conditions with $\delta(G) \ge \frac{3}{4}n$ admits a $K_3$-decomposition. In an addendum to Nash-Williams' article~\cite{nash1970unsolved}, Graham remarked that the fraction $3/4$ would be best possible. After decades of partial progress, the Nash-Williams' conjecture was very recently proven true in a breakthrough result by the the first and fourth author~\cite{DP-NW}. As for $K_q$-decompositions of sparse graphs, in 2019 Glock, K\"{u}hn, Lo, Montgomery, and Osthus~\cite{GKLMO19} showed that the existence of a \emph{fractional} $K_q$-decomposition is enough to guarantee the existence of a $K_q$-decomposition, and in 2019 Montgomery~\cite{montgomery2019fractional} produced the best-known result on the fractional relaxation by obtaining a minimum degree bound of $(1-\frac{1}{100q})n$. In 2025, the authors~\cite{DHLP25Counter} disproved a folklore generalization of the Nash-Williams' conjecture by showing that, for each integer $q\geq 4$, there exist a $c>1$ and infinitely many $K_q$-divisible graphs $G$ with minimum degree at least $\big(1-\frac{1}{c\cdot (q+1)}\big)v(G)$ with no $K_q$-decomposition and even no fractional $K_q$-decomposition, thus leaving a wide gap between known lower and upper bounds for the minimum degree $K_q$-decomposition threshold.

As for clique decompositions of sparse hypergraphs, in 2021, Glock, K{\"u}hn, and Osthus~\cite{GKO20Survey} posited that every large enough $r$-graph $G$ that satisfies the necessary divisibility conditions with $\delta_{r-1}(G) \geq \left(1 - \Theta_{r}\left(\frac{1}{q^{r-1}}\right)\right)\cdot v(G)$ admits a $K_q^r$-decomposition. Recent work by the authors improves on work from Glock, K{\"u}hn, Lo, and Osthus \cite{GKLO16} by bounding this $K_q^r$-decomposition threshold on almost the same order as posited by the ``Hypergraph Nash-Williams' conjecture.'' In particular, the first, third, and fourth authors~\cite{DLP25} showed that, for every $r \geq 3$ and real $\varepsilon\in (0,1]$, there exists a constant $c > 0$ such that minimum codegree at least $1-\frac{c}{q^{r-1+\varepsilon}}$ suffices for the existence of a \emph{fractional} $K_q^r$-decomposition in large enough $r$-graphs. Meanwhile, the second and fourth authors~\cite{HP25} recently tied the $K_q^r$-decomposition threshold to its fractional relaxation, thus implying that $\delta_{r-1}(G) \ge \left(1 -\frac{c}{q^{r-1+\eps}}\right) \cdot n$ suffices for the existence of a $K_q^r$-decomposition in a large enough $r$-graph $G$ satisfying the necessary divisibility conditions.

For the $K_q^r$-subgraph setting of Tur\'an theory, the $K_q^r$-factor setting of a hypergraph Hajnal-Szemer\'edi theorem, and the $K_q^r$-decompositions discussed above, 
there exists an $i \in \{0, 1, r\}$ such that the respective result provides a density threshold to guarantee there exists a collection of $K_q^r$-subgraphs of a divisible graph $G$ such that every $i$-set of $G$ is in exactly one element in the set (where we view $G$ as a simplicial complex). Since a copy of $K^r_q$ contains the empty set, Tur\'an's theorem answers the $i = 0$ case, while the hypergraph Nash-Williams' conjecture addresses the $i = r$ case. Thus,~\cref{thm:mainexact} serves as a stepping stone between the influential areas of Tur\'an theory and combinatorial designs by addressing this natural simplicial complex question in the case where $i = 1$. 

\subsection{Barriers to Hypergraph Clique Factors}\label{Subsec:Barriers}
We utilize a result of Lang as a black box to streamline our proof of~\cref{thm:mainexact}.
%or example, Lo and Markstr\"{o}m~\cite{LM15} used absorption for their 2015 upper bond on the $K_q^r$-factor minimum codegree threshold, Lu and Sz\'ekely~\cite{LS07} employ the Lov\'asz Local Lemma in the context of random injections to prove~\cref{thm:PreviousHHS}, and Mycroft~\cite{M11} utilized hypergraph regularity and a hypergraph blow-up lemma for his work on 
In 2023, Lang~\cite{lang2023tiling} showed that the asymptotic minimum degree threshold for the existence of a hypergraph factor is governed by three obstructions: space, divisibility, and covering. Though Lang's results extend to $F$-factors for any hypergraph $F$, we include the simplified definitions and statements restricted to the case where $F= K_q^r$. 

A \textit{fractional $K_q^r$-factor} of an $r$-graph $G$ is an assignment of weights $\omega: \mc{K}_q^r(G) \rightarrow [0, 1]$ such that for every vertex $v\in V(G)$, we have $\sum_{K \in \mc{K}_q^r(G): v\in V(K)}\omega(K) = 1$, where $\mathcal{K}_q^r(G)$ denotes the set of subgraphs of $G$ isomorphic to $K_q^r$.
Likewise, an \textit{integral $K_q^r$-factor} of $G$ is an assignment of weights $\omega: \mc{K}_q^r(G) \rightarrow \mathbb{Z}$ such that for every vertex $v\in V(G)$, we have $\sum_{K \in \mc{K}_q^r(G): v\in V(K)}\omega(K) = 1$. Finally, a \textit{$K_q^r$-covering} of $G$ is an assignment of weights $\omega: \mc{K}_q^r(G) \rightarrow \{1\}$ such that for every vertex $v\in V(G)$, we have $\sum_{K \in \mc{K}_q^r(G): v\in V(K)}\omega(K) \ge 1$ (in other words, we have that every vertex is contained in at least one copy of $K_q^r$). Note that the existence of a fractional factor, an integral factor, and a covering are each trivially necessary for the existence of a factor in $G$. 

For ease of discussion, we introduce the following notation. For all $1 \le r < q$, let $\delta_{q,r}$ be the infimum over all $\delta\in[0,1]$ such that for $n$ large enough with $q \mid n$, every $n$-vertex $r$-graph $G$ with $\delta_{r-1}(G) \ge \delta \cdot  n$ admits a $K_q^r$-factor. 
(We note that in~\cite{lang2023tiling}, this threshold is defined such that $\delta_{r-1}(G) \ge (\delta_{q, r} +\eps)n$ guarantees the existence of a factor, but omitting $\eps$ provides an equivalent definition given the use of ``infimum.'') 
Analogously define $\delta^{Spa}_{q,r}$ for fractional $K_q^r$-factors, $\delta^{Div}_{q,r}$ for integral $K_q^r$-factors, and $\delta^{Cov}_{q,r}$ for $K_q^r$-coverings.
Lang showed that $\delta_{q,r}$ is tied to these space, divisibility, and covering thresholds.

\begin{thm}[Lang~\cite{lang2023tiling}]
	For integers $1\leq r < q$, we have
	\[
		\delta_{q,r} = \max \big\{\delta^{Spa}_{q,r},	\,\delta^{Div}_{q,r},\, \delta^{Cov}_{q,r} \big\}.
	\]
\end{thm}

This result is simplified in our hypergraph clique factor setting. In particular, the existence of a fractional $K_q^r$-factor (where every vertex receives a total weight of one) implies that each vertex is contained in at least one copy of $K_q^r$, and hence $\delta^{Cov}_{q,r}\leq\delta^{Spa}_{q,r}.$ Thus, in our setting, Lang's result simplifies to:
%Explain why this is true for cliques but not nec in general?
\begin{cor}
    For integers $1\leq r < q$, we have
	\[
		\delta_{q,r} = \max \big\{\delta^{Spa}_{q,r},	\,\delta^{Div}_{q,r}\big\}.
	\]
\end{cor}

In this paper, we bound $\delta^{Div}_{q,r}$ and $\delta^{Spa}_{q,r}$ as follows.

\begin{lem}\label{lem:DivBarrier}
 For every $r\geq 3$, there exists a constant $c>0$ such that for every integer $q>r$ we have
    \[\delta^{Div}_{q,r} \le 1-c \cdot {\frac{\log q}{q^{r-1}}}.\]
\end{lem}

\cref{lem:DivBarrier} is a straightforward argument achieved by bounding $\delta^{Div}_{q,r}$ by a related connectivity threshold which is, in turn, satisfied in our high minimum degree setting. Thus, bounding the fractional threshold is key to determining $\delta_{q,r}$ and, hence, the main result of this paper is the following. 

\begin{thm}\label{thm:main}
For every $r\geq 2$, there exists a constant $c>0$ such that for every integer $q>r$ we have
\[\delta^{Spa}_{q,r}\leq 
1-c\cdot\frac{(\log q)^{1-\frac{1}{r-1}}}{q^{r-1}}.\] 
\end{thm}

\cref{thm:main} is a significant improvement on the fractional threshold for $K_q^r$-factors, as the previous best known upper bound appears to have been given by~\cref{thm:PreviousHHS} (which also holds for fractional factors, as shown by~\cref{cor:FracHHR}).
\cref{thm:mainexact} now follows as a corollary of~\cref{lem:DivBarrier} and~\cref{thm:main} (we note that the $r = 2$ case follows from~\cref{thm:PreviousHHS}). We detail the novel ideas in our proof of~\cref{thm:main} in~\cref{Subsec:Novelty} and provide an overview of the argument in~\cref{Sec:PfOverview}. We also note that in our proof we require a technical theorem which we can repeatedly apply to find a large independent set in a sparse hypergraph setting. We introduce such a result, \cref{thm:maincoloring}, in~\cref{subsec:indepsets}.

\subsection{On the Novelty of Our Proof of the Fractional Threshold}\label{Subsec:Novelty}

To prove~\cref{thm:main} we utilize a linear programming framework. Inspired by the first and fourth author's~\cite{DP-NW} proof of the Nash-Williams' conjecture, we show that if an $r$-graph $G$ does not admit a fractional $K_q^r$-factor, then $G$ admits a vertex weighting that satisfies a non-negative local condition (the sum of weight of the vertices in any given copy of $K_q^r$ is non-negative) and a negative global condition (the sum of the weight of every vertex in $V(G)$ is negative). Our argument is novel because of how we exploit this local condition to contradict the global one in our hypergraph setting by playing the non-negative vertices against the negative ones. We direct readers to~\cref{Sec:PfOverview} for a more detailed proof overview of~\cref{thm:main}.

Utilizing a linear programming perspective in this context is natural and indeed has been used for past results in similar settings. For instance, R\"{o}dl, Ruci\'{n}ski, and Szemer\'edi's~\cite{RRS09} 2009 hypergraph matching threshold, Keevash and Mycroft's~\cite{KM15Geometric} 2015 work on a geometric theory for hypergraph matchings, and Bowtell, Kathapurkar, Morrison, and Mycroft's~\cite{BKMM25} 2025 result for generalized triangle factors all use Farkas' Lemma, a classical linear programming result, to reach a contradiction via direct counting arguments. In fact, one can even use a dual perspective to obtain a short proof of the fractional version of the Hajnal-Szemer\'edi theorem: by translating the statement of a fractional version of~\cref{thm:HS} into a linear program, we get that its dual minimizes the sum of non-negative vertex weights on $V(G)$ such that the weights of vertices in any given copy of $K_q$ sum to at least one. Starting with the most negatively weighted vertex $v_1$, we iteratively build a $q$-clique $K$ on vertices $v_1, v_2, \ldots, v_q$ where each $v_i$ is a lowest weight common neighbor of $v_1, \ldots, v_{i-1}.$ Then, after ordering the vertices in non-decreasing weight, the minimum degree condition on $G$ implies that each $v_i$ appears as the first vertex of the $i$th interval or earlier. Since the sum of the weight in $K$ is at least one, the total weight of $V(G)$ is at least $n/q$ and the desired fractional $K_q$-factor follows from strong duality. 
This argument is a special case of Keevash and Mycroft's Lemma 3.6 in~\cite{KM15Geometric}. Generalizing this proof to the $r$-graph setting is not enough to prove~\cref{thm:main} (it instead gives a codegree threshold on the order of $(1 - 1/q^{r-1})n$), but we are nevertheless inspired to consider a set $S$ of the most negatively weighted vertices and employ discharging techniques to spread positive weight to $S$, thus lower bounding the weight in $V(G)$. 

Our first innovation is in how we handle the negative vertices. We leverage codegree Tur\'an theory to show that all of the negative weight is concentrated in a small proportion of vertices. It therefore stands to reason there are many cliques of non-negative vertices that extend to a copy of $K_q^r$ with a given clique of negative ones. Yet, we still cannot guarantee that the most negative vertices form a clique. We sidestep this obstacle by taking a fractional $K_p^r$-factor of the negative vertices for an appropriately chosen value $p<q$. This implies the existence of a $p$-clique $T$ of negatives whose average weight is at least (in absolute value) the average of all negatives, thus enabling us to lower bound the negative weight without necessarily taking the lowest-weighted vertices.

Our next ingenuity is how we show the non-negative vertices outweigh the negatives. 
First, given the copy $T$ of $K_p^r$ of negatively weighted vertices of average weight, we only want to spread weight to $T$ from copies of $K_q^r$ that otherwise consist of only non-negative vertices. Inspired by the second and fourth author's work with rooted hypergraph Tur\'an theory in~\cite{HP25}, we consider a non-uniform, auxiliary hypergraph $H$ on only the non-negatively weighted vertices where an edge of uniformity $i$ represents a set of $i$ vertices that does not extend to an edge in $G$ with some $(r-i)$-subset of $V(T)$. Here, an independent set of size $q-p$ represents the desired clique in $G$. 
To show the non-negative weight of $G$ outweighs the negative, we only want to count the weight of each non-negative vertex at most once. Therefore, we find enough disjoint independent sets in $H$, each of size $q-p$, using~\cref{thm:maincoloring}, our original independent set result for non-uniform sparse hypergraphs. We detail the history and novelty of~\cref{thm:maincoloring} in~\cref{sec:indep}. This concludes the proof, provided the number of disjoint sets is large enough relative to the number of negative vertices over $p$, and we choose $p$ such that this holds.

% The crux of our proof of~\cref{thm:main} is not just how we individually solve these challenges but how we overcome them simultaneously. In particular, through our careful choice of the parameter $p$ we are able to balance the competing needs for a fractional $K_p^r$-decomposition within the negative vertices and the existence of many disjoint, large sets of size $(q - p)$ that extend to a clique within the positive vertices to obtain the desired contradiction. 

\subsection{Independent Sets in Sparse Hypergraphs}\label{subsec:indepsets}

The principal technical result of the paper is an independent-set theorem for sparse non-uniform hypergraphs. It is used in two distinct parts of the argument. In~\cref{sec:divisibility}, it produces common clique extensions
of pairs of vertices, allowing us to control the divisibility barrier. In the
proof of~\cref{thm:main}, it is applied iteratively to produce many disjoint
extensions of a fixed clique, which is the key step in bounding the fractional threshold. We believe that it may be of independent interest as a non-uniform independent-set result for sparse hypergraphs.

For a hypergraph $H$, an \emph{independent set} is a set of vertices containing
no edge of $H$, and the independence number $\alpha(H)$ is the maximum size of
an independent set in $H$. A proper coloring of $H$ is a partition of $V(H)$
into independent sets, and the chromatic number $\chi(H)$ is the minimum number
of parts in such a partition. Consequently, $\alpha(H)\cdot \chi(H)\geq v(H)$, and every upper bound on $\chi(H)$ immediately yields a corresponding lower bound on $\alpha(H)$.
Recall that the maximum $i$-degree of a hypergraph $G$ is the maximum, over all vertices $v$ of $G$, of the number of edges of $G$ with size $i$ containing $v$, that is 
\[\Delta_i(G):=\max_{v\in V(G)} |\set{e\in G\colon |e|=i,\ v\in e}.\]

A standard application of the Lov\'{a}sz Local Lemma (see e.g. Cooper and Mubayi~\cite[Lemma~8]{CM16}) yields that every $r$-bounded hypergraph $G$ with maximum $i$-degree $\Delta_i$, $i=2,\ldots,r$ satisfies \(\chi(G)=O_r\parens*{\max_{i=2}^r\ \Delta_i^{1/(i-1)}}\).
There is a very rich history of improvements of this general bound under additional sparsity assumptions. Ajtai, Koml{\'o}s, Pintz, Spencer, and Szemer{\'e}di~\cite{AKPSS82} obtained a logarithmic improvement (ie. $\log^{1 / (r-1)}\Delta_{r}$) for $r$-uniform hypergraphs containing no short Berge cycles. (The proof of this statement goes through for non-uniform hypergraphs, and modern concentration techniques may facilitate the application of this result in lower girth hypergraphs. Even still, it is still not enough for the low codegree setting of~\cref{thm:mainexact}.) Duke, Lefmann, and R{\"o}dl~\cite{DLR95} established the same order of magnitude for linear hypergraphs. More generally, Kostochka, Mubayi, and Verstra{\"e}te~\cite{KMV14} derived lower bounds on the independence number from bounds on the maximum codegree. 
Further generalizations of these results to hypergraph colorings were obtained 
by Frieze and Mubayi~\cite{FM2013} for linear hypergraphs, by Cooper and Mubayi~\cite{CM2015} for rank-three hypergraphs, and more recently by Li and the fourth author~\cite{LP22} for hypergraphs of arbitrary bounded rank. These results illustrate that substantially stronger conclusions become possible when the edges of the hypergraph have limited overlap. We therefore introduce the following notions of codegrees. For an $r$-bounded hypergraph $H$ and $k\in[r]$, let $H^{(k)}$ be the $k$-uniform subhypergraph of $H$ containing all $k$-edges of $H$. We then define the \emph{maximum $(k,\ell)$-codegree} of $H$ to be 
\[\Delta_{\ell}(H^{(k)}) = \max_{S\in \binom{V(H)}{\ell}} |\{e\in E(H): S\subseteq e, |e|=k\}|.\]
Similarly, the \emph{common $2$-degree} of distinct vertices $u, v\in V(H)$ is defined as $|\{w\in V(H): uw, vw\in E(H)\}|$. The \emph{maximum common $2$-degree} of $H$ is the maximum of the common $2$-degree of $u$ and $v$ over all distinct pairs of vertices $u,v$ of $H$.

The following theorem is the closest existing result to the setting needed in this article. It is a corollary of the aforementioned result of Li and the fourth author.
\begin{thm}[Li and Postle~\cite{LP22}]\label{thm:DP}
Let $r\ge 3$ be an integer and $\beta > 0$ be real. There exists $D_{\beta}$ such that the following holds for all $D\ge D_{\beta}$:

Let $H$ be a $r$-bounded hypergraph whose edges have size at least two such that $\Delta_i(H) \le O(D^{i-1}\log D)$ for all $2\le i \le r$. If the maximum common $2$-degree of $H$ is at most $D^{1-\beta}$ and for all $2\le \ell < k \le r$, we have that
$\Delta_{\ell}(H^{(k)}) \le D^{k-\ell-\beta}$,
then 
$$\chi(H) \leq D.$$
\end{thm}
We note that~\cref{thm:DP} also follows from the first and fourth author's~\cite{DP22} 2022 work on hypergraph matchings avoiding forbidden submatchings. 
It is also implied by more recent work on conflict-free matchings, such as that of Glock, Joos, Kim, K\"
{u}hn, and Lichev~\cite{GJKKL22}.

This result almost works for our setting, except that the uniformity two edges pose a critical bottleneck. In particular, as described in~\cref{Subsec:Novelty}, in order to find rooted cliques in our $r$-graph $G$, we want to find independent sets in a non-uniform, auxiliary setting. Therefore, though we only apply this argument when $r \ge 3$, our auxiliary hypergraph includes 2-uniform edges and we have no control over the maximum common 2-degree. Thus, we prove the following result, which unites the standard, greedy graph coloring bound and the above low codegrees bound. 

\begin{thm}\label{thm:maincoloring}
    For every integer $r\geq 3$ and reals $\beta,K > 0$, there exist $c=c(r,\beta)\in(0,1)$ and $q_0=q_0(r,\beta,K)$ such that for all $q\geq q_0$ there exists $D_0=D_0(q)$ such that the following holds for all $D\geq D_0$. 

    Let $H$ be an $r$-bounded hypergraph on $N$ vertices with $N\leq KqD$, and whose edges have size at least two. If for all integers $2\leq \ell<j\leq r$ we have
    \[\begin{cases}
        \Delta_1(H^{(2)})&\leq D\\
        \Delta_{1}(H^{(j)})&\leq D^{j-1}\cdot\log q\\
        \Delta_{\ell}(H^{(j)})&\leq D^{j-\ell}\cdot q^{-\beta}
    \end{cases}\]
    then we have
    \[\alpha(H)\geq c\cdot \frac{N}{D}.\]
\end{thm}

Though this statement is sufficient for our purposes, it is possible that the large $D$ condition is not necessary in the statement of~\cref{thm:maincoloring}. Additionally, we conjecture that the coloring and list coloring generalizations of this result hold true. 

To prove~\cref{thm:maincoloring}, we avoid rerunning the proof of~\cref{thm:DP} as we do not care about exact bounds. Instead, we use conflict free machinery from Glock, Joos, Kim, K\"{u}hn, and Lichev~\cite{GJKKL22} with test functions to track which 2-uniform edges survive after random sparsification before greedily choosing a sufficiently large independent set. We quarantine the technicalities of this result to~\cref{sec:indep}, where we include both a more in-depth overview of the argument as well as the proof itself. 

\subsection{Notation and Paper Organization}
The remainder of the paper is organized as follows. We provide a proof overview for~\cref{thm:main} in~\cref{Sec:PfOverview}. In~\cref{sec:divisibility}, we focus on the divisibility barrier for $K_q^r$-factors and prove~\cref{lem:DivBarrier}. In~\cref{sec:space}, we turn our attention to our main result about fractional thresholds by proving~\cref{thm:main}. Finally, in~\cref{sec:indep}, we prove~\cref{thm:maincoloring}, thus establishing our technical independent set result.

Throughout the remainder of the paper, we require some additional notation. Let $G$ be a (hyper)graph. If $S \subseteq V(G)$, then $G - S$ denotes the subgraph of $G$ induced by $V(G) \setminus S$. Next, $\binom{V(G)}{i}:=\{S\subseteq V(G) \; : \; |S| = i \}$. We use $H \subseteq G$ to denote that $H$ is a subgraph of $G$ and $\mathcal{K}_q^r(G)$ denotes the set of subgraphs of $G$ isomorphic to $K_q^r$. For a vertex weighting $\phi:V(G) \rightarrow \mathbb{R}$ and a subset $U \subset V(G)$, let $\phi(U):= \sum_{v \in U}\phi(v)$. Finally, $[n]_i:= \{i, ..., n\}.$ 

\section{\texorpdfstring{\cref{thm:main}}{Theorem~\ref{thm:main}} Proof Overview}\label{Sec:PfOverview}
As discussed in~\cref{Subsec:Novelty}, we use a linear programming framework in our proof of~\cref{thm:main} introduced by the first and fourth authors in~\cite{DP-NW} for their proof of the Nash-Williams' conjecture. In line with the vocabulary introduced in~\cite{DP-NW}, given a (hyper)graph $G$, an assignment of weights $\phi:V(G) \rightarrow \mathbb{R}$ is \textit{$K_q^r$-covered} if for all $H\in \mc{K}_q^r(G),$ we have that $\sum_{v\in V(H)}\phi(v)\geq 0$. Likewise, we say $\phi$ is \textit{totally-negative} if $\sum_{v\in V(G)} \phi(v) < 0$. For a vertex weighting $\phi$ of $V(G)$, we use $V_{\phi}^-(G)$ and $V_{\phi}^+(G)$ to denote the negatively weighted and non-negatively weighted vertices of $G$, respectively. When $G$ and $\phi$ are implicit, we write $V^-$ and $V^+$ for brevity.
The following is an immediate extension of the classical result of Farkas' Lemma~\cite{FarkasLem}, as stated by the first and fourth authors in~\cite{DP-NW} in the context of graph decompositions, and adapted for $K_q^r$-factors here. 

\begin{lem}[Farkas' Lemma for $K_q^r$-factors]\label{lem:NoFactor}
    Let $G$ be an $r$-graph on $n$ vertices. Then $G$ admits no fractional $K_q^r$-factor if and only if there exists an assignment of weights $\phi:V(G)\to [-1,+\infty)$ such that $\phi$ is $K_q^r$-covered and totally-negative.
\end{lem}
\cref{lem:NoFactor} follows directly from Farkas' Lemma applied to the vertex-clique incident matrix of an $r$-graph. As discussed in~\cref{Subsec:Novelty}, we exploit this local non-negativity condition to contradict the global negativity condition of~\cref{lem:NoFactor} to prove~\cref{thm:main}. To do so, we spread weight from $V^+$ to $V^-$ by first lower bounding the negative weight in $G$ and then showing that the non-negatively weighted vertices can afford to spread their weight to the negative ones. To that end, we consider an $r$-graph $G$ on $n$ vertices with $\delta_{r-1}(G)\geq \parens*{1-c\cdot \gamma}\cdot n$ for $\gamma:= \frac{\log^{1-1/(r-1)}q}{q^{r-1}}$ and assume for a contradiction that $G$ does not admit a fractional $K_q^r$-factor. 

First, in order to lower bound the negative weight in $G$, we show there is a clique of average weight in $G[V^-]$. To that end, for a vertex weighting $\phi$ of $G$ and a subset $S \subseteq V(G)$ we use $\overline{\phi}(S) : = \frac{1}{|S|} \sum_{v \in S}\phi(v)$. We prove that, for an appropriate value $p$, there exists a copy $T$ of $K_p^r$ in $G[V^-]$ with $|\overline{\phi}(T)| \ge |\overline{\phi}(V^-)|$, which we then use as a representative set of negative vertices in our subsequent arguments. The question of which value of $p$ to use is thus as the heart of our argument, and our choice balances two competing arguments. On the one hand, we use the existence of a fractional $K_p^r$-factor to find the copy $K_p^r$ of average weight. Therefore, $p$ should be small enough so that the missing codegree inherited by $G[V^-]$ implies the existence of a fractional $K_p^r$-factor, as guaranteed by~\cref{cor:FracHHR}, a corollary of the previous codegree $K_q^r$-factor threshold of Lu and Sz\'ekely~\cite{LS07}. On the other hand, $p$ should be as large as possible: a copy $T$ of $K_p^r$ in $G[V^-]$ with $|\overline{\phi}(T)| \ge |\overline{\phi}(V^-)|$ carries at least a $p/|V^-|$-proportion of the total negative weight in $G$. Hence, the larger $p$ is, the more leverage we have in the final discharging argument. Therefore, we take $p$ essentially maximal subject to
\[  c\cdot \gamma \cdot n \leq \frac{\sqrt{c}}{p^{r-1}}\cdot |V^-|.\]
%We note that if $|V^-| \ge (r+1)^{r-1} \sqrt{c} \cdot \gamma \cdot n,$ then our choice of $p$ ensures $p \ge r+1$. If $|V^-|$ is small, enough negative weight is concentrated in the single lowest weight vertex in $V^-$. If $|V^-|$ is large, the split in the definition below ensures that we always obtain $p\geq r+1$.
We note that we use a ceiling function to ensure we take $p \ge 1$. However, it is possible that we choose $p < r$, in which case a copy of $K_p^r$ is an independent set on $p$ vertices.

\begin{lem}\label{lem:averageKpr}
For every $r\geq 3$, there exists a constant $c>0$ such that for every integer $q>r$ the following holds for every large enough $n$. Let $\gamma:=q^{1-r} \cdot \log^{1-1/(r-1)}q$, $G$ be an $r$-graph on $n$ vertices with $\delta_{r-1}(G) \ge \left(1 -c\cdot\gamma\right)n$, and $\phi:V(G) \rightarrow \mathbb{R}$ be a vertex weighting of $G$ such that $V^-$ is nonempty. For
\begin{align*}
        p:= \floor*{\parens*{\ceil*{\frac{|V^-|}{\sqrt{c} \cdot \gamma \cdot n}}}^{1/(r-1)}}
\end{align*}
% \begin{align*}
    %     p:= 
    %     \begin{cases}
    %     \floor*{\parens*{\frac{|V^-| }{\sqrt{c} \cdot \gamma \cdot n}}^{1/(r-1)}}
    %     &\text{ if } |V^-|\geq (r+1)^{r-1} \sqrt{c} \cdot \gamma \cdot n,\\
    %     1&\text{ otherwise}.
    %     \end{cases}
    %\end{align*}
there exists a copy $T$ of $K_p^r$ in $G[V^-]$ such that $|\overline{\phi}(T)| \ge | \overline{\phi}(V^-)|$.
\end{lem}
To prove~\cref{lem:averageKpr}, we utilize a corollary of~\cref{thm:PreviousHHS}, the previous codegree $K_q^r$-factor threshold of Lu and Sz\'ekely~\cite{LS07}, to show that $G[V^-]$ admits a fractional $K_p^r$-factor which we then use to find our desired copy of $K_p^r$. 
We state and prove this corollary, \cref{cor:FracHHR}, and prove~\cref{lem:averageKpr} in~\cref{Subsec:PfRFracDecomp}.

The other key piece of our proof of~\cref{thm:main} is a result about finding many disjoint copies of $K_{q-p}^r$ in $G[V^+]$ that extend to a $q$-clique in $G$ with a given copy $T$ of $K_p^r$ in $G[V^-].$ In this way, we prevent any non-negative vertex from being overburdened when we spread weight from $V^+$ to $T$.
In this way, the weight of the non-negative vertices outweighs the negative weight of $T$, and hence $V^-$, yielding a contradiction. To that end, for a given copy $T$ of $K_p^r$, we say a set of vertices $S \subseteq V(G)$ is a \emph{$K_q^r$-extension} of $T$ if $S$ induces a copy of $K_{q-p}^r$ and $S \cup T$ induces a copy of $K_q^r$. Additionally, we say two $K_q^r$-extensions $S_1$ and $S_2$ are \emph{disjoint} if $S_1 \cap S_2 = \emptyset$. 

\begin{lem}\label{lem:IndepSets}
    For every $r\geq 3$, there exists a constant $c>0$ such that for every large enough integer $q>r$ the following holds for every large enough $n$. Let $\gamma:=q^{1-r} \cdot \log^{1-1/(r-1)}q$, $G$ be an $r$-graph on $n$ vertices with  $\delta_{r-1}(G) \ge \left(1 -c \cdot \gamma\right)n$, and $\phi:V(G) \rightarrow [-1, +\infty)$ be a $K_q^r$-covered vertex weighting of $G$ such that $V^-$ is nonempty. Let 
    \begin{align*}
        p:= \floor*{\parens*{\ceil*{\frac{|V^-|}{\sqrt{c} \cdot \gamma \cdot n}}}^{1/(r-1)}}
\end{align*}
    If $T$ is a copy of $K_p^r$ in $G[V^-],$ then there exists a family $\mathcal{F}$ of at least $|V^+|/(4q)$ disjoint $K_q^r$-extensions of $T$ in $V^+$.
\end{lem}
Proving~\cref{lem:IndepSets} is the main challenge of our argument. As described in~\cref{Subsec:Novelty}, we accomplish this by applying~\cref{thm:maincoloring}, our technical result about large independent sets in sparse non-uniform hypergraphs, to an auxiliary hypergraph where independent sets represent cliques in $G[V^+]$ that extend with $T.$ 
In order to show that we can iteratively apply~\cref{thm:maincoloring} to obtain enough of these independent sets, we also leverage codegree Tur\'an theory to show that the negative weight in $V(G)$ is concentrated in a small number of vertices (and hence $|V^+|$ is large) as follows.

\begin{lem}\label{Lem:RSmall}
For every $r\geq 3$, there exists a constant $c>0$ such that for every integer $q>r$ the following holds for every large enough $n$. Let $\gamma:=q^{1-r} \cdot \log^{1-1/(r-1)}q$. If $G$ is an $r$-graph on $n$ vertices with  $\delta_{r-1}(G) \ge \left(1 -c\cdot \gamma\right)n$ and $\phi:V(G) \rightarrow [-1, +\infty]$ is a $K_q^r$-covered vertex weighting of $G$, then 
\[|V_{\phi}^-(G)| <  \frac{\sqrt{c} \cdot n}{\log^{1/(r-1)}{q}}.\]
\end{lem}
\begin{proof}
Let $c \in (0, 1)$ be such that~\cref{thm:HyperTuran} holds for uniformity $r$ with parameter $\sqrt{c}$. 
We assume $|V^-| \ge c  \cdot n/\log^{1/(r-1)}{q}$, as otherwise the result trivially holds since $c \in (0, 1)$. As $\phi$ is $K_q^r$-covered, we have that $G[V^-]$ is $K_q^r$-free. Since $n$ is large enough, by~\cref{thm:HyperTuran},
    \[\delta_{r-1}(G[V^-])<\left( 1-\sqrt{c}\cdot\frac{\log q}{q^{r-1}}\right)|V^-|,\]
  while the minimum $(r-1)$-degree of $G$ implies
  \[\delta_{r-1}(G[V^-])\geq |V^-|-c\cdot \gamma \cdot n.\]
Combining these two inequalities, we obtain $|V_{\phi}^-(G)| <  \sqrt{c} \cdot n/\log^{1/(r-1)}{q}$, as desired.
\end{proof}
Given~\cref{Lem:RSmall}, we prove~\cref{lem:IndepSets} in~\cref{Subsec:PfIndepSets}.

Assuming the above lemmas, the proof of~\cref{thm:main} now follows by showing that we can indeed use the disjoint clique extensions from~\cref{lem:IndepSets} to spread weight from $V^+$ to a copy $T$ of $K_p^r$ in $G[V^-]$ of average weight guaranteed by~\cref{lem:averageKpr}, lower bounding the weight of all of $V^-$ and thus contradicting the fact that $\phi$ is totally-negative.
We include the proof of~\cref{thm:main} in~\cref{Subsec:MainPf}.

\section{Divisibility Barrier}\label{sec:divisibility}

Recall that $\delta^{Div}_{q,r}$ is the minimum $(r-1)$-degree threshold that ensures that for every large enough $r$-graph $G$ satisfying $q \mid v(G)$ and $\delta_{r-1}(G)\geq \delta^{Div}_{q,r}  \cdot n,$ we have that $G$ admits an integral $K_q^r$-factor. In this section, we prove~\cref{lem:DivBarrier} by relating $\delta^{Div}_{q,r}$ to the notion of the connectivity threshold of $K_q^r$. 

\subsection{Connectivity Threshold}

In~\cite{lang2023tiling}, Lang introduced a notion of hypergraph connectivity which he used to bound $\delta^{Div}_{q,r}.$ For an $r$-graph $G$ and $\ell \geq r$, let $K_{\ell}(G)$ denote the (2-uniform) graph on $V(G)$ obtained by placing a $2$-uniform $\ell$-clique on the vertex set of every $r$-uniform $\ell$-clique of $G$. The hypergraph $G$ is \textit{$\ell$-connected} if $K_{\ell}(G)$ is connected (in the usual graph theoretic sense). Let $\delta_{q, r}^{Con}$ denote the \textit{connectivity threshold of $K_q^r$}, that is, the minimum $(r-1)$-degree threshold that forces $q$-connectivity in an $r$-uniform host graph. Lang related this parameter to $\delta^{Div}_{q,r}$ as follows.

\begin{lem}[Special case of Corollary 5.5 in~\cite{lang2023tiling}]\label{lem:DivThreshLeqConThresh}
    Let $q>r\geq3$. Then
\[\delta^{Div}_{q,r}\leq \delta^{Con}_{q,r}.\]
\end{lem}
In fact,~\cref{lem:DivThreshLeqConThresh} is a special case of Lang's stronger result concerning all hypergraphs satisfying certain properties, which includes $K_q^r$. 
To show~\cref{lem:DivBarrier}, it therefore suffices to show that $\delta^{Con}_{q,r} \leq 1-c \cdot {\frac{\log q}{q^{r-1}}}$ for some constant $c$ depending only on $r$. To that end, we include the following definition. 

\begin{definition}\label{def:extendable}
    Let $q \ge 2$ be an integer. We say a graph is \textit{$q$-extendable} if for every $u, v \in V(G)$ there exists a copy $K$ of $K_{q-1}^r$ that extends with both $u$ and $v$, that is, $K \cup \{u\} \cong K_{q}^r$ and $K \cup \{v\} \cong K_{q}^r$.
\end{definition}
 % \CH{(We could insert a parameter for number of vertices and make it $(2, q)$-extendable if we want.)}

 Note that if a hypergraph is $q$-extendable, then it is trivially $q$-connected. Thus, in order to bound $\delta^{Con}_{q,r}$ (and, hence, $\delta^{Div}_{q,r}$) we show that our high minimum degree setting satisfies this stronger condition by proving the following lemma.
 
\begin{lem}\label{lem:CliqueExtension}
For every integer $r \ge 3$, there exists a real $C > 0$ such that for all integers $q > r$, the following holds for every large enough $n$. If $G$ is an $r$-graph on $n$ vertices with $\delta_{r-1}(G) \ge \left(1-C\cdot\frac{\log q}{q^{r-1} }\right)n$, then $G$ is $q$-extendable. 
\end{lem}

To prove~\cref{lem:CliqueExtension}, we fix vertices $u, v \in V(G)$ and work in an auxiliary, non-uniform hypergraph consisting of non-edges of $G$ and sets of $(r-1)$ vertices that do not extend to an edge with both $u$ and $v$. We then apply~\cref{thm:maincoloring} to find an independent set of size at least $q-1$ in this auxiliary graph, which corresponds to the desired clique in $G$.

\begin{proof}
Fix integers $3 \le r < q$. Let $n$ be large enough. Fix $\beta:=1/2$ and let $c\in(0,1)$ be a real number and $q_0$ be an integer such that~\cref{thm:maincoloring} holds for every integer $q' \ge q_0$ when $K = 1/c$. 

We claim it is enough to prove the result for $q \ge q_0$. If $q < q_0$, assume \cref{lem:CliqueExtension} holds for $q_0$ with constant $C_{q_0}$. In this case, let $C := C_{q_0} \cdot \frac{\log q_0}{q_0^{r-1}}\cdot \frac{r^{r-1}}{\log r}.$ Note that $C > 0$ and, since $q_0$ is a constant depending on $r$, $C$ is also a constant depending only on $r$. Since $r < q$, $C = C_{q_0}\cdot \frac{\log q_0}{q_0^{r-1}}\cdot \frac{r^{r-1}}{\log r} \le C_{q_0}\cdot \frac{\log q_0}{q_0^{r-1}}\cdot \frac{q^{r-1}}{\log q}.$ This implies $C\cdot\frac{\log q}{q^{r-1} } \le C_{q_0}\cdot\frac{\log q_0}{q_0^{r-1} }$, and hence $\delta_{r-1}(G) \ge \left(1 - C_{q_0}\cdot\frac{\log q_0}{q_0^{r-1}} \right)n$. In this case, it thus suffices to apply the argument for $q_0$: any copy of $K_{q_0 - 1}^r$ that extends with two vertices $u$ and $v$ contains the desired $K_{q-1}^r$ subgraph.

Therefore, assume $q \ge q_0$. Let $D_0$ be as in \cref{thm:maincoloring}. Let $D:= c(n-2)/(q-1)$, and observe that $D \ge D_0$ since $n$ is large enough. Finally, let $C:= c^{r-1}$ (and note that $c < 1$ implies $C \le c^{i}$ for all $0 \le i \le r-1$). Now let $G$ be an $r$-bounded hypergraph with $\delta_{r-1}(G) \ge \left(1-C\cdot\frac{\log q}{q^{r-1} } \right)n$. Consider distinct vertices $u$ and $v$ in $V(G)$.

Now we define the auxiliary, non-uniform hypergraph $H$ to which we will apply~\cref{thm:maincoloring}. To that end, let 
\[N_u := \{ S \subseteq V(G) : |S| = r-1 \text{ and } S \cup \{u\} \notin E(G) \}.\]
Similarly define $N_v$. 
Let $H$ be the non-uniform hypergraph with $V(H) = V(G) \setminus \{u, v\}$ and $E(H) = E(\overline{G}) \cup N_u \cup N_v$. Observe that it suffices to find an independent set of size at least $q-1$ in $H$, as this corresponds to a copy of $K_{q-1}^r$ that extends with both $u$ and $v$ in $G$. 

Note that $H$ is an $r$-bounded hypergraph whose edges have size at least $r-1 \ge 2$. Also, observe that $v(H) = n-2 \le Kqc(n-2)/(q-1) = KqD$ by the choice of $K$. Therefore, in order to employ \cref{thm:maincoloring} to find an appropriately large independent set in $H$, it suffices to check the relevant degree conditions. We do so by considering the ``missing'' edges of $G$.
To that end, let $w \in V(G)$. There are at most $\binom{n}{r-2}$ subsets of $V(G)$ of size $(r-1)$ containing $w$, and, by the minimum degree of $G$, at most $C \cdot \log q \cdot n/q^{r-1}$ edges of $\overline{G}$ containing such an $(r-1)$-set. Thus, by the choice of $D$ and $C$ and the fact that $n$ is large enough,
\[\Delta_1(H^{(r)}) \le \frac{\binom{n}{r-2} \cdot C \cdot  \log q \cdot n}{q^{r-1}} \le  \frac{Cn^{r-1} \log q}{q^{r-1}} \le D^{r-1}\log q. \]
Likewise, for $x \in \{u, v\}$, there are at most $\binom{n}{r-3}$ subsets of $V(G)$ of size $(r-1)$ containing $w$ and $x$. Again by the minimum degree of $G$, choice of $D$, and choice of $C$, we get that
\[\Delta_1(H^{(r-1)}) \le \frac{C n^{r-2}\log q}{q^{r-1}} \le D^{r-2}\log q.\]
In particular, if $r-1 = 2$, this means
\[\Delta_1(H^{(2)}) = \Delta_1(H^{(r-1)}) \le \frac{Cn^{r-2}\log q}{q^{r-1}} = \frac{Cn\log q}{q^{2}} \le \frac{c (n-2)}{q-1} = D,\]
where the last inequality is because $C \le c$, $\log q/q^2 \le 1/(q-1)$, and $n$ is large enough.
If instead $r-1 \ge 3$, then $\Delta_1(H^{(2)}) = 0 \le D.$

For the final two calculations, note that since $\log q/q < 1$ and $1/q^{r-2} < 1/(q-1)^{r - 2.5}$, we have that $\log q/ q^{r-1} < 1/(q-1)^{r - 2.5}$. Thus, for each integer $i$ with $0 \le i \le r-3$, since $C \le c^{i+1}$ and $n$ is large enough, we have that 
\begin{align}
    \frac{C \cdot n^{i+1}\log q}{ q^{r-1}} \le \frac{c^{i+1} (n-2)^{i+1}}{(q-1)^{r - 2.5}} \le \frac{c^{i+1} (n-2)^{i+1}}{(q-1)^{i+3/2}} \le D^{i+1} \cdot q^{-\beta}. \label{Observation}
\end{align}

Now, let $\ell := r-1 -i$ where $i$ is an integer with $0 \le i \le r-3$. 
There are at most $\binom{n}{i}$ subsets of $V(G)$ of size $(r-1)$ containing a given set of $\ell$ vertices. Again by the minimum degree of $G$ and choice of $D$, we have
\[\Delta_{\ell}(H^{(r)})\le \frac{Cn^{i+1} \log q}{q^{r-1}}  \le D^{i+1} \cdot q^{-\beta} = D^{r-\ell} \cdot q^{-\beta},\]
where the last inequality is a consequence of \cref{Observation}. Note that if $r = 3$ we are done.
Similarly, if $r \ge 4$, let $\ell' := r - 2 - i'$ where $i'$ is an integer with $0 \le i' \le r-4.$ For each choice of $x \in \{u, v\}$, there are at most $\binom{n}{i}$ subsets of $V(G)$ of size $(r-1)$ containing $x$ and a given set of $\ell$ additional vertices. Once more by the minimum degree of $G$ and choice of $D$, we have
\[\Delta_{\ell'}(H^{(r-1)}) \le \frac{Cn^{i'+1} \log q}{q^{r-1}}  \le D^{i'+1} \cdot q^{-\beta} = D^{r - 1-\ell'} \cdot q^{-\beta},\]
where the last inequality is a consequence of \cref{Observation}

Therefore, by \cref{thm:maincoloring}, we have that $\alpha(H) \ge c \cdot v(H)/D = c(n-2)/D = q-1$ by our choice of $D$, as desired.
\end{proof}

Since~\cref{lem:CliqueExtension} implies that $\delta^{Con}_{q,r} \leq 1-c \cdot {\frac{\log q}{q^{r-1}}}$ for a constant $c$ depending only on $r$, the proof of~\cref{lem:DivBarrier} now follows from~\cref{lem:DivThreshLeqConThresh}.

\section{Space Barrier - Improving the Fractional Threshold}\label{sec:space}

In this section, we focus our attention on~\cref{thm:main}, our improved upper bound on fractional $K_q^r$-factors. First, we prove~\cref{lem:averageKpr} in~\cref{Subsec:PfRFracDecomp} and~\cref{lem:IndepSets} in~\cref{Subsec:PfIndepSets} before combining these into a proof of~\cref{thm:main} in~\cref{Subsec:MainPf} using the methodology described in~\cref{Sec:PfOverview}.

\subsection{Proof of
~\texorpdfstring{\cref{lem:averageKpr}}{Lemma~\ref{lem:averageKpr}}
}\label{Subsec:tools}\label{Subsec:PfRFracDecomp}
Our proof of~\cref{lem:averageKpr} utilizes the following fractional version of Lu and Sz\'ekely~\cite{LS07} codegree $K_q^r$-factor threshold result, \cref{thm:PreviousHHS}. \cref{cor:FracHHR} is an immediate corollary of~\cref{thm:PreviousHHS}. We include its proof here for completeness. 

\begin{cor}\label{cor:FracHHR}
    For every integer $r\geq 2$, there exists a constant $c > 0$ such that the following holds for every $q>r$. Every large enough $r$-graph $G$ with 
    \[\delta_{r-1}(G) \geq \left(1 - \frac{c}{q^{r-1}}\right)\cdot v(G),\]
    contains a fractional $K_q^r$-factor.
\end{cor}

\begin{proof}
Let $q>r\geq 2$, let $n$ be large enough with $n\equiv i\pmod q$  for some $i\in[q-1]_0$. Let $c_0>0$ such that~\cref{thm:PreviousHHS} holds, and let $c=c_0/2$. Let $G$ be a $n$-vertex
$r$-graph with \(\delta_{r-1}(G) \geq \big(1 - \frac{c}{q^{r-1}}\big)\cdot n\). Fix a set $S$ of $i$ vertices in $V(G)$. Observe that
\[\delta_{r-1}(G-S)\geq \delta_{r-1}(G)-|S|\geq \left(1 - \frac{c}{q^{r-1}}\right)\cdot n-q\geq \left(1 - \frac{c_0}{q^{r-1}}\right)\cdot (n-|S|),\]
as $n$ is large enough. Therefore, by~\cref{thm:PreviousHHS}, let $\varphi_S$ be a $K_q^r$-factor of $G-S$.
% that is, an assignment of nonnegative weights to every clique $H\in\mathcal{K}_q^r(G-S)$ such that for every vertex $u$ in $V(G)\setminus S$, we have
% \[\sum_{H\in \mathcal{K}_q^r(G-S)\colon u\in H} \varphi_S(H)=1.\]
For simplicity, let $\varphi_S(H)=0$ for every clique $H\in\mathcal{K}_q^r(G)\setminus\mathcal{K}_q^r(G-S)$. Then, for every $H\in\mathcal{K}_q^r(G) $, let 
\[\varphi(H):=\frac{1}{\binom{n-1}{i}}\cdot\sum_{S\in \binom{V(G)}{i}} \varphi_S(H).\]
For every vertex $u\in V(G)$, we obtain 
\begin{align*}
\sum_{H\in \mathcal{K}_q^r(G)\colon u\in H} \varphi(H)
&=\frac{1}{\binom{n-1}{i}}\cdot \sum_{H\in \mathcal{K}_q^r(G)\colon u\in H} \sum_{S\in \binom{V(G)}{i}} \varphi_S(H)\\
&=\frac{1}{\binom{n-1}{i}}\cdot \sum_{S\in \binom{V(G)}{i}} \sum_{H\in \mathcal{K}_q^r(G)\colon u\in H} \varphi_S(H)\\
&=\frac{1}{\binom{n-1}{i}}\cdot \sum_{S\in \binom{V(G)}{i}} \mathbf{1}_{u\not\in S}\\
&=1,
\end{align*}
as desired.
\end{proof}

With~\cref{cor:FracHHR} in hand, we are now prepared to prove~\cref{lem:averageKpr}.

\begin{lateproof}{lem:averageKpr}
   % Let $c \in (0, 1)$ be such that~\cref{cor:FracHHR} holds with parameter $\sqrt{c}$. We first claim that $G[V^-]$ has a fractional $K_p^r$-factor. If $p < r$, this is trivially true, \CH{What if $p = r.$ Need to cite an additional source here I think.} Now assume $|V^-|\geq (r+1)^{r-1} \sqrt{c}\cdot \frac{\log^{1-1/(r-1)}q}{q^{r-1}} \cdot n.$
   % Using the minimum codegree of $G$ and the definition of $p$, we deduce that 
   %  \[\delta_{r-1}(G[V^-])\geq|V^-|-c\cdot \frac{\log^{1-1/(r-1)}q}{q^{r-1}} \cdot n 
   %  \geq \Big(1-\frac{\sqrt{c}}{p^{r-1}}\Big)|V^-|.\]
   %  Since $n$, and hence $|V^-|$, is large enough,~\cref{cor:FracHHR} implies that $G[V^-]$ admits a fractional $K_p^r$-factor $\psi$.

    Let $a_r$ be the constant from~\cref{cor:FracHHR}, choose $c\in(0,1)$ so that $2^r\sqrt{c}\leq a_r$, and let $x:=|V^-|/(\sqrt{c}\gamma n)$. Since $|V^-|\geq1$ and $\sqrt{c}\gamma n\geq1$ for sufficiently large $n$, we have $1\leq p\leq\ceil*{x}\leq |V^-|$.
    
    If $p<r$, we give each $p$-subset of $V^-$ a weight of $\binom{|V^-|-1}{p-1}^{-1}$. As each vertex lies in $\binom{|V^-|-1}{p-1}$ such sets, this yields a fractional $K_p^r$-factor of $V^-$. If $p\geq r$, set $s:=\max\{p,r+1\}$. Since $p^{r-1}\leq\ceil*{x}<x+1$ and $p\geq r$, we have $x>1$ and $s^{r-1}\leq(2p)^{r-1}\leq2^r x$. Thus
    \[ \delta_{r-1}(G[V^-])\geq |V^-|-c\gamma n =\parens*{1-\frac{\sqrt{c}}{x}}|V^-|
    \geq\parens*{1-\frac{a_r}{s^{r-1}}}|V^-|.\]
    For fixed $r,q,c$, the integer $s$ is bounded and $|V^-|>(r^{r-1}-1)\sqrt{c}\gamma n$. Therefore~\cref{cor:FracHHR} yields a fractional $K_s^r$-factor in $G[V^-]$. If $p>r$, then $s=p$. If $p=r$, let each $(r+1)$-clique of weight $w$ contribute $w/r$ to each of its $r+1$ edges. Each vertex of the clique lies in exactly $r$ of these edges and receives total weight $w$ from the clique. In all cases, we obtain a fractional $K_p^r$-factor $\psi$ of $V^-$.

    Let $\mc{P}$ be the set of all $K_p^r$ subgraphs in $G[V^-]$. Since $\psi$ is a fractional $K_p^r$-factor, we have
    \[\sum_{T\in\mc{P}}\psi(T) =\frac1p\sum_{v\in V^-} \sum_{\substack{T\in\mc{P}\\v\in V(T)}}\psi(T)
    =\frac{|V^-|}{p}.\]
    Similarly,
    \[ \sum_{T\in\mc{P}}\psi(T)\overline{\phi}(T)
    =\frac1p\sum_{v\in V^-}\phi(v) \sum_{\substack{T\in\mc{P}\\v\in V(T)}}\psi(T)
    =\frac{|V^-|}{p} \overline{\phi}(V^-).\]
    Hence some $T\in\mc{P}$ with $\psi(T)>0$ satisfies $\overline{\phi}(T)\leq\overline{\phi}(V^-)<0$, and therefore $|\overline{\phi}(T)|\geq|\overline{\phi}(V^-)|$, as desired.
\end{lateproof}

\subsection{Proof of~
\texorpdfstring{\cref{lem:IndepSets}}{Lemma~\ref{lem:IndepSets}}
}\label{Subsec:PfIndepSets}
Now we prove~\cref{lem:IndepSets} as follows. 

\begin{lateproof}{lem:IndepSets}
    Let $r\geq 3$ be an integer, let $\beta:=\frac12$, and let $c=c(r)\in(0,1)$ be an arbitrarily small constant such that~\cref{Lem:RSmall} holds with parameter $c$, \cref{thm:maincoloring} holds with parameter $\sqrt{c}$, and $64c\leq 1$.
    Let $K:=2/\sqrt{c}$, and let $q_0=q_0(r)$ be an integer such that~\cref{thm:maincoloring} holds. Let $q$ be large enough compared to $1/c$, hence compared to $r$, in particular such that $q\geq q_0$. Finally, let $n$ be a sufficiently large integer, and, for brevity of calculations, define
    \[L:=\log^{1/(r-1)}q,\qquad \gamma:= \frac{\log^{1- 1/(r-1)}q}{q^{r-1}} = \frac{\log q}{L\cdot q^{r-1}}.\]
    Observe that
    \begin{equation}\label{eq:epsL}
        L^{r-2}=\gamma \cdot q^{r-1}.
    \end{equation}

    Put $x:=|V^-|/(\sqrt{c}\gamma n)$ and $B:=q/L$. By~\cref{Lem:RSmall}, $x<1/(\gamma L)=B^{r-1}$. Thus, for sufficiently large~$q$,
    \begin{equation}\label{eq:Boundp}
        p^{r-1}\leq\ceil*{x}<B^{r-1}+1\leq\tfrac43 B^{r-1}.
    \end{equation}
    
 % Combining the definition of $p$ and the bound on $|V^-|$ obtained by~\cref{Lem:RSmall} we get that $p^{r-1}<\frac{1}{\gamma L} = \frac{q^{r-1}}{L^{r-1}}$ and hence, as $L<q$,
 %    \begin{equation}\label{eq:Boundp}
 %    p<\frac{q}{L}.
 %    \end{equation}

    % Suppose first that $p=1$. We therefore have $|V^-|< (r+1)^{r-1} c\cdot \gamma \cdot n$, and, for $q$ large enough,
    % \[N:=|V^+|\geq (1-(r+1)^{r-1} c \cdot \gamma) n \geq \frac{n}{2}.\]

    Suppose first that $p=1$. By~\cref{Lem:RSmall}, for sufficiently large~$q$ we have \[ N:=|V^+|>n(1-\sqrt{c}/L)\geq n/2.\]
        
    Recall~\cref{def:extendable} and let $C$ be such that~\cref{lem:CliqueExtension} holds. Let $T=\{v\}$ for some vertex $v\in V^-$, and let $V^+_0:=V^+$. We construct the required family greedily. Suppose that for some integer $i\geq0$, the set $V^+_i\subseteq V^+$ satisfies $|V^+_i|\geq N/2$. Let $G_i:=G[V^+_i\cup\{v\}]$. We aim to apply~\cref{lem:CliqueExtension} to $G_i$. To that end, observe that $v(G_i)\geq \frac{N}{2}\geq \frac{n}{4}$. 
    Therefore we obtain 
    \[ \delta_{r-1}(G_i)\geq v(G_i)-c \cdot \gamma \cdot n \geq (1-4c\gamma)v(G_i).\]
    For $q$ large enough we have
    %\[L:=\log^{1/(r-1)}q,\qquad \gamma:=\frac{\log q}{L\cdot q^{r-1}}.\]
    \[4c\gamma=4c\cdot \frac{\log^{1-1/(r-1)} q}{q^{r-1}} 
    < C\cdot\frac{\log q}{q^{r-1} }.\]
    Therefore, by~\cref{lem:CliqueExtension}, $G_i$ is $q$-extendable. By~\cref{def:extendable}, there exists a set $S_i\subseteq V^+_i$, with $|S_i|=q-1$, and such that $G[S_i\cup\{v\}]\cong K_q^{r}$. Let $V^+_{i+1}:=V^+_i\setminus S_i$, and continue until fewer than $N/2$ vertices remain. We obtain a family of $m$ pairwise disjoint sets with $m(q-1)>\frac{N}{2}$, hence $m>\frac{N}{2(q-1)}\geq \frac{N}{4q}$, as desired.

    Assume now that $p>1$. Let $V^+_T\subseteq V^+$ be the set of vertices in $V^+$ such that, for all $u\in V^+_T$, we have that $T\cup \{u\}$ induces a $(p+1)$-clique in $G$, that is,
    \[V^+_T:=\set*{u\in V^+ \; \colon \; G[T\cup \{u\}]\cong K_{p+1}^r},\]
    and let $N:=|V^+_T|$.\bigskip
    
    % Recall that $T$ induces a $p$-clique in $G$, and, by the minimum degree of $G$, that any set of $(r-1)$ vertices in $T$ does not extend to an edge with at most \(c\cdot \gamma n\) vertices in $G$. Therefore, using the fact that $p = \floor*{\parens*{\frac{|V^-|}{\sqrt{c}\gamma n}}^{1/(r-1)}}$, the fact that $|V^-|< \frac{\sqrt{c}\cdot n}{L}$ by~\cref{Lem:RSmall}, and the definition of $c$, we have
    % \[|V^+_T|\geq |V^+|-\binom{p}{r-1}\cdot c\cdot \gamma  n
    % \geq n-|V^-|-\frac{1}{2}\cdot p^{r-1}\cdot c\cdot \gamma  n
    % \geq n-|V^-|-\frac{\sqrt{c}\cdot |V^-|}{2}
    % \geq n\cdot\parens*{1-\frac{\sqrt{c}}{L}-\frac{c}{2L}}\geq \frac{n}{2}.\]

    Each $(r-1)$-set in $T$ has at most $c\gamma n$ non-neighbors. Since $\binom{p}{r-1}\leq p^{r-1}$ and $\gamma B^{r-1}=1/L$,~\cref{Lem:RSmall,eq:Boundp} yield,
    \[N\geq n-|V^-|-\binom{p}{r-1}c\gamma n
    \geq n\parens*{1-\frac{\sqrt{c}}{L}-\frac{4c}{3L}}\geq\frac{2n}{3}.\]    

    Let $\mc{H}$ be the $r$-bounded hypergraph with $V(\mc{H})=V_T^+$ and
    \[E(\mc{H}):=\set*{S\subseteq V^+_T\colon |S|\in\{2,\ldots,r\},~\exists \; T'\in \binom{T}{r-|S|}\text{ s.t. }T'\cup S\not\in E(G)},\]
    that is, a $k$-edge $U$ of $\mc{H}$ is a $k$-set of non-negatively weighted vertices in $V^+_T$, i.e. every individual vertex of $U$ induces a copy of $K_{p+1}^{r}$ with $T$ but together $U\cup T$ do not form a clique in $G$. In particular, we have that for a set $S \subseteq V^+$, if $S$ is independent in $\mc{H}$ then $G[T\cup S]$ is complete. 
    
    Fix integers $1\leq \ell<j\leq r$, and let $X\subseteq V^+_T$ be an $\ell$-set. We now compute $\Delta_{\ell}(\mc{H}^{(j)})$. To do so, recall that by the definition of $\mc{H},$ for every $j$-edge $S$ containing $X$, there exists $T'\in\binom{T}{r-j}$ such that $T'\cup S\not\in E(G)$. We then select $j-\ell-1$ other vertices in $V^+_T$ to form an $r-1$-set. By the minimum degree assumption of $G$, that $(r-1)$-set has at most $c\gamma n$ non-neighbors in $G$ (and hence at most that many neighbors in $\mc{H}$). Therefore
    \[\Delta_{\ell}(\mc{H}^{(j)})\leq \binom{|T|}{r-j}\cdot \binom{|V^+_T|}{j-\ell-1}\cdot c\gamma n.\]

    Since $n\leq3N/2$ and $|T|=p$, we obtain
    \[\Delta_{\ell}(\mc{H}^{(j)}) \leq\frac32\cdot p^{r-j}N^{j-\ell}c\gamma.\]
    By~\cref{eq:Boundp}, $p<q$ for sufficiently large $q$, and
    $p^{r-j}\leq(4/3)B^{r-j}$ for $2\leq j\leq r$.
    Using~\cref{eq:epsL}, we obtain
    % As $2|V^+_T|\geq n$, $|T|=p$, and $|V^+_T|=N$, we obtain
    % \[\Delta_{\ell}(\mc{H}^{(j)})\leq 2p^{r-j}\cdot N^{j-\ell}\cdot c\gamma.\]        
    % Recall that, by~\cref{eq:epsL} we have $L^{r-2}=\gamma q^{r-1}$, and that $pL\leq q$ by~\cref{eq:Boundp}, therefore we obtain
    \begin{align*}
    \Delta_{\ell}(\mc{H}^{(j)})
    &\leq 2\cdot\parens*{\frac{q}{L}}^{r-j}\cdot N^{j-\ell}\cdot c\cdot \frac{L^{r-2}}{q^{r-1}}
    \leq 2c\cdot\frac{1}{q^{j-1}}\cdot L^{j-2}\cdot N^{j-\ell}.        
    \end{align*}
    Recall that $\beta=\frac12$ and $K=\frac{2}{\sqrt{c}}$, and let
    \[D:=\frac{N}{Kq}=\frac{\sqrt{c}N}{2q}.\]
    For $j=2$ and $\ell=1$ we obtain 
    \[ \Delta_{1}(\mc{H}^{(2)}) \leq 
    2c\cdot\frac{1}{q}\cdot N
    \leq \frac{\sqrt{c}N}{2q}
    =D,\]
    where we used $64c\leq 1$, thus $4c\leq\sqrt{c}$.
    For $j\in\{3,\ldots,r\}$ and $\ell=1$ we obtain, for $q$ large enough compared to $1/c$,
    \[\Delta_{1}(\mc{H}^{(j)})
    \leq 2c\cdot\frac{1}{q^{j-1}}\cdot L^{j-2}\cdot N^{j-1}
    \le 2^j\cdot\frac{1}{(\sqrt{c})^{j-3}}\cdot \parens*{\frac{\sqrt{c}N}{2q}}^{j-1}\cdot(\log q)^{(r-2)/(r-1)}
    \leq D^{j-1}\cdot\log q.\]
    For $2\leq\ell<j\leq r$ we obtain, for $q$ large enough compared to $1/c$,
    \[\Delta_{\ell}(\mc{H}^{(j)}) \leq 
    2c\cdot\frac{1}{q^{j-1}}\cdot L^{j-2}\cdot N^{j-\ell}
    \leq 2^j\cdot\frac{1}{(\sqrt{c})^{j-\ell-2}}\cdot \parens*{\frac{\sqrt{c}N}{2q}}^{j-\ell}\cdot \frac{\log q}{q}
    \leq D^{j-\ell}\cdot q^{-\beta}.
    \]
    To summarize, we have obtained that for all integers $2\leq \ell<j\leq r$ we have
    \[\begin{cases}
        \Delta_1(\mc{H}^{(2)})&\leq D\\
        \Delta_{1}(\mc{H}^{(j)})&\leq D^{j-1}\cdot\log q\\
        \Delta_{\ell}(\mc{H}^{(j)})&\leq D^{j-\ell}\cdot q^{-\beta},
    \end{cases}\]
    which are exactly the requirements of~\cref{thm:maincoloring}. Let $D_0(q)$ be given by~\cref{thm:maincoloring}. As $n$ is a large enough integer compared to $q$, and using $2N\geq n$ and $D=N/Kq$, we obtain that $D\geq D_0(q)$. Therefore, we aim to apply~\cref{thm:maincoloring} iteratively at least $m:=\floor*{N/(2q)}$ times. Let $V^+_0:=V^+_T$. Let $i\in[m]_0$, assume that $V^+_i\subseteq V^+_T$ satisfies $|V^+_i|\geq N/2$, and let $\mc{H}_i:=\mc{H}[V^+_i]$. We trivially have, for all $1\leq \ell<j\leq r$,
    \[\Delta_{\ell}(\mc{H}_i^{(j)})\leq \Delta_{\ell}(\mc{H}^{(j)}),\]
    and $|V^+_i|\leq |V^+_T|=KqD$. Therefore we can apply~\cref{thm:maincoloring} to $\mc{H}_i$ and we obtain
    \[\alpha(\mc{H}_i)\geq \sqrt{c}\cdot\frac{|V^+_i|}{D}\geq \sqrt{c}\cdot \frac{N}{2D}=q.\]
    Let $S_i\subseteq V^+_i$ be an independent set of size $q-p$ in $\mc{H}_i$, and let $V^+_{i+1}:=V^+_i\setminus S_i$. Observe that $i\leq m$ ensures that 
    \[|V^+_i|\geq |V^+_T|-qm =  N - q\cdot \floor*{\frac{N}{2q}}\geq \frac{N}{2}.\]
    Let $\mc{F}=\set*{S_0,\ldots,S_{m}}$, with $m+1=\floor*{N/(2q)}+1\geq N/(2q)\geq n/(4q) \ge |V^+|/4q$. Each $S_i$ is an independent set in $\mc{H}$, therefore $S_i\cup T$ induces a copy of $K_q^r$ in $G$, as desired.
\end{lateproof}

\subsection{Proof of Main Fractional Result}\label{Subsec:MainPf}
Now we are prepared to prove~\cref{thm:main} using the method outlined in~\cref{Sec:PfOverview}.

\begin{lateproof}{thm:main}
We note that the statement is true for $r=2$, as guaranteed by~\cref{thm:PreviousHHS}. We further note that it suffices to prove the statement for $q$ sufficiently large compared to $r$. Indeed, assume that the result holds for every $q'>q_0$, for some fixed $q_0=q_0(r)$. Then, by decreasing the constant $c$ if necessary, for all $q\in\{r+1,\ldots,q_0\}$ we have 
    \[c\cdot \log^{1-1/(r-1)}q\leq c\cdot \log^{1-1/(r-1)}q_0\leq c',\]
    where $c'$ is obtained from~\cref{cor:FracHHR}. Then we obtain that, for all $q\in\{r+1,\ldots,q_0\}$, we have
    \[1-c\cdot \frac{\log^{1-1/(r-1)} q}{q^{r-1}}\geq 1-\frac{c'}{q^{r-1}},\]
    and the result in that range of $q$ also follows immediately from~\cref{cor:FracHHR}.

    Thus, let $r \ge 3$ and let $c \in(0, 1)$ be an arbitrarily small constant such that~\cref{lem:averageKpr}, \cref{lem:IndepSets}, and~\cref{Lem:RSmall} all hold with parameter $c$, $4\sqrt{c}\leq \log^{1/(r-1)}3$, and $64c\leq 1$. Let $q > r$ be sufficiently large, and 
let $G$ be an $r$-graph on $n$ vertices with  $\delta_{r-1}(G) \ge \left(1 -c\cdot \frac{\log^{1-1/(r-1)}q}{q^{r-1}}\right)n$. Assume for a contradiction that $G$ does not admit a fractional $K_q^r$-factor. Then, by~\cref{lem:NoFactor}, there exists an assignment of weights to the vertices of $G$, $\phi:V(G) \rightarrow [-1, +\infty)$, such that $\phi$ is $K_q^r$-covered and totally-negative.
% Let 
% \begin{equation}\label{eq:def_p}
%         p:= 
%         \begin{cases}
%         \floor*{\parens*{\frac{|V^-| \cdot q^{r-1}}{\sqrt{c}\cdot \log^{1-1/(r-1)}q \cdot n}}^{1/(r-1)}}
%         &\text{ if } |V^-|\geq (r+1)^{r-1} \sqrt{c} \cdot \frac{\log^{1-1/(r-1)}q}{q^{r-1}} \cdot n,\\
%         1&\text{ otherwise}.
%         \end{cases}
%     \end{equation}
Let $\gamma:=q^{1-r}\log^{1-1/(r-1)}q$ and set
\begin{equation}\label{eq:def_p}
    x:=\frac{|V^-|}{\sqrt{c}\gamma n},\qquad
    p:=\floor*{\parens*{\ceil*{x}}^{1/(r-1)}}.
\end{equation}
Since $\phi$ is totally-negative, $V^-$ is nonempty. Therefore, by~\cref{lem:averageKpr} there exists a copy $T$ of $K_p^r$ in $G[V^-]$ such that $|\overline{\phi}(T)| \ge | \overline{\phi}(V^-)|$. By~\cref{lem:IndepSets}, there exists a family $\mathcal{F}$ of at least $|V^+|/(4q)$ disjoint subsets of $V^+$, each of size $q-p$, and such that for every $S\in\mathcal{F}$, we have that $S\cup T$ induces a copy of $K_q^r$ in $G$. Since $\phi$ is $K_q^r$-covered, for every $S \in \mc{F}$ we have that $\phi(S) \ge |\phi(T)|$. Thus, 
    \begin{equation}\label{eq:SpreadWeight}
        \phi(V^+)\geq \sum_{S \in \mc{F}}\phi(S) \ge \frac{|V^+|}{4q}\cdot|\phi(T)| = \frac{n -|V^-|}{4q}\cdot|\phi(T)|.
    \end{equation}

Therefore, we have that
\[|\phi(V^-)| = |\overline{\phi}(V^-)| \cdot |V^-|  
     \le |\overline{\phi}(T)| \cdot |V^-| 
     = \frac{|\phi(T)|}{|T|} \cdot |V^-|
     \le \frac{4q}{ n-|V^-|} \cdot \frac{|V^-|}{p} \cdot \phi(V^+),\]
% \begin{align*}
%     |\phi(V^-)| &= |\overline{\phi}(V^-)| \cdot |V^-| \\ 
%     & \le |\overline{\phi}(T)| \cdot |V^-| \\
%     & = \frac{|\phi(T)|}{|T|} \cdot |V^-| \\
%     & \le \frac{4q}{|V^+|} \cdot \frac{|V^-|}{p} \cdot \phi(V^+),
% \end{align*}
where the first inequality follows from the choice of $T$ and the second inequality follows from~\cref{eq:SpreadWeight}.  It now suffices to show that $\frac{p}{4q} \cdot \frac{n-|V^-|}{|V^-|} \ge 1$, as this implies $|\phi(V^-)| \le \phi(V^+)$, contradicting the fact that $\phi$ is totally negative. 

% With this in mind, first suppose $p = 1$. By~\cref{eq:def_p}, we have $|V^-| < (r+1)^{r-1} \sqrt{c} \cdot  \log^{1 - 1/(r-1)}q \cdot n/q^{r-1} \le  n/2$ as $q$ is large enough. Using this and the fact that $ r \ge 3$, we get 
With this in mind, first suppose $p=1$. By~\cref{eq:def_p}, $\ceil*{x}<2^{r-1}$, so $|V^-|<2^{r-1}\sqrt{c}\gamma n\leq n/2$ for sufficiently large $q$. Using $r\geq3$, we get
\[\frac{p}{4q} \cdot \frac{n - |V^-|}{|V^-|} \ge \frac{n}{8q \cdot |V^-|} \ge \frac{q^{r-2}}{8\cdot 2^{r-1}\sqrt{c} \cdot \log^{1 - 1/(r-1)}q} \ge 1,\]
where the last inequality is because $q$ is large enough compared to $r$ (and hence also $c$). 

Now assume $p > 1$. Recall that by the choice of $c$, 
\begin{equation}\label{eq:cL}
4\sqrt{c}\leq \log^{1/(r-1)}q.
\end{equation}
Since $p\geq2$,~\cref{eq:def_p} yields $x^{1/(r-1)}\leq\ceil*{x}^{1/(r-1)}<p+1\leq3p/2$. By~\cref{eq:cL}, $1-\sqrt{c}/\log^{1/(r-1)}q\geq3/4$. Hence $p(1-\sqrt{c}/\log^{1/(r-1)}q)>x^{1/(r-1)}/2$. Therefore,
\begin{align*}   
       \frac{p}{4q}\cdot\frac{n-|V^-|}{|V^-|} & >\frac{p}{4q}\cdot\frac{n}{|V^-|}\cdot\parens*{1-\frac{\sqrt{c}}{\log^{1/(r-1)}q}}&\qquad \text{ by \cref{Lem:RSmall}}\\  
       &>\frac{1}{8q}\cdot \parens*{\frac{|V^-|\cdot q^{r-1}}{\sqrt{c}\cdot \log^{1 - 1/(r-1)}q \cdot  n}}^{1/(r-1)}
       \cdot\frac{n}{|V^-|}&\qquad \text{ by }\cref{eq:cL}\text{ and }\cref{eq:def_p}\\  
       &=\frac{1}{8q}\cdot \parens*{\frac{n}{|V^-|}}^{1-1/(r-1)}
       \cdot \parens*{\frac{q^{r-1}}{\sqrt{c}\cdot \log^{1 - 1/(r-1)}q}}^{1/(r-1)}\\
       &>\frac{1}{8q}\cdot \parens*{\frac{\log^{(r-2)/(r-1)}q  \cdot q^{r-1}}{\log^{1-1/(r-1)}q}}^{1/(r-1)}\cdot\frac{1}{\sqrt{c}}&\qquad \text{ by \cref{Lem:RSmall}}\\
       &=\frac{1}{8}\cdot\frac{1}{\sqrt{c}}\\   
       &\geq 1 &\qquad \text{ as }64c\leq1,
    \end{align*}
as desired. Hence, $|\phi(V^-)| \le \phi(V^+)$, contradicting~\cref{lem:NoFactor} and thus $G$ admits a fractional $K_q^r$-factor. 
\end{lateproof}

\section{Independent sets in non-uniform hypergraphs}\label{sec:indep}

In this section, we turn our attention to the proof of~\cref{thm:maincoloring}, where we prove that a hypergraph $H$ with good control over its codegrees contains a large independent set. We begin with a proof overview in~\cref{Subsec:ColorPfOverview}, before including the relevant conflict free definitions and statements in~\cref{subsec:ConflictFree} and finally our proof of~\cref{thm:maincoloring} in~\cref{subsec:IndepPf}.

\subsection{\texorpdfstring{\cref{thm:maincoloring}}{Theorem~\ref{thm:maincoloring}} Proof Overview}\label{Subsec:ColorPfOverview}

As discussed in~\cref{subsec:indepsets}, 
we require a result about large independent sets in sparse, non-uniform hypergraphs in our proof of~\cref{lem:IndepSets}. However, we have no control over the maximum common 2-degree of our non-uniform auxiliary graph, and so none of the existing machinery for independent sets in hypergraphs is equipped to handle the 2-uniform edges in our proof of~\cref{lem:IndepSets}. 
Surprisingly, \cref{thm:maincoloring} combines the greedy bound for uniformity two edges while maintaining the existing sparse bounds for larger uniformities given by Li and Postle's~\cref{thm:DP}. Since we do not optimize the constant in~\cref{thm:maincoloring}, instead of rerunning the technicalities of the proof of~\cref{thm:DP} we streamline the proof of~\cref{thm:maincoloring} by utilizing conflict-free matching results with test functions as follows. 

To prove~\cref{thm:maincoloring}, we first pass to a random induced subhypergraph $H'\subseteq H$ with $v(H')=\Theta(q^2N/D)$ and still well controlled (co)degrees. It then suffices to find an independent set in $H'$ of order $\Omega(v(H')/q^2)$. We encode partial list-colorings of $H'$ as matchings in an auxiliary graph. We then use the conflict machinery with test-systems of Glock, Joos, Kim, K\"
{u}hn, and Lichev~\cite{GJKKL22}. Informally, their theorem produces a large matching that avoids a forbidden family of submatchings, while retaining the expected number of submatchings from a specified family (called a test system). In our setting, this provides a partial coloring of $H'$ with few monochromatic edges, where we use test-functions to keep track of the 2-uniform edges. Here we encounter an annoying technicality: pairs of vertices with too many common neighbors pose a bottleneck to finding the desired independent set. We overcome this by adding these pairs to the conflict system. 
Finally, after some careful alterations to delete the remaining monochromatic edges, we obtain a partial proper coloring of $H'$, and the largest color class yields the required independent set.

\subsection{Conflict Free Machinery}\label{subsec:ConflictFree}

In this section, we include the definitions and statements necessary for our proof of~\cref{thm:maincoloring}.

First, in order to formally state Glock, Joos, Kim, K\"
{u}hn, and Lichev's~\cite{GJKKL22} conflict-free theorem, we require the following definitions and notation. For a vertex $v\in V(\mc{H})$, we denote by $N_{\mc{H}}^{(2)}(v)$ the {\em neighborhood} of $v$ in $\mc{H}$, that is, 
\[N_{\mc{H}}^{(2)}(v):=\set*{u\in V(\mc{H})\colon \{u,v\}\in E(\mc{H})},\]
and by $\mc{H}_v$ the {\em link} of $v$ in $\mc{H}$, that is, the hypergraph with 
\[V(\mc{H}_v)=V(\mc{H})\setminus\{v\},\quad E(\mc{H}_v)=\set*{e\setminus\{v\}\colon e\in\mc{H},\ v\in e}.\]
For an integer~$j\geq1$, recall that we write $\mc{H}^{(j)}$ for the $j$-uniform subgraph of~$\mc{H}$ containing all edges of size exactly $j$. For an integer~$j\geq 1$ and~$v\in V(\mc{H})$, we write~$\mc{H}^{(j)}_v$ as a shorthand for~$(\mc{H}_v)^{(j)}=(\mc{H}^{(j+1)})_v$.

For a graph $G$, a {\em conflict system for $G$} is a (not necessarily uniform) hypergraph~$\mc{C}$ with vertex set $E(G)$ whose edges have size at least $2$. We call the edges of $\mc{C}$ {\em conflicts} of $\mc{C}$ and for~$e\in G$ and~$C\in\mc{C}_e$, we say that~$C$ is a {\em semiconflict}. We say that an edge-set $E\subseteq E(G)$ is $\mc{C}$\emph{-free} if no edge of $\mc{C}$ is a subset of $E$. To construct a $\mc{C}$-free matching of $G$, we require control over $\mc{C}$ in the following sense:

\begin{definition}[bounded conflict-systems]\label{def:BoundedConflict}
    Let $G$ be a graph, and $\mc{C}$ be a conflict system for $G$. For integers~$d\geq 1$ and~$\ell\geq 2$ and reals~$\Gamma\geq 0$ and~$\eps\in(0,1)$, we say that~$\mc{C}$ is $(d,\ell,\Gamma,\eps)$-{\em bounded} if all the following properties hold.
	\begin{enumerate}[label=(C\arabic*)]
		\item $2\leq\abs{C}\leq \ell$ for all~$C\in \mc{C}$;\label{item: conflict size} 

		\item $\sum_{j\in[\ell]} \frac{\Delta_{1}(\mc{C}^{(j)})}{d^{j-1}}\leq \Gamma$ 
        and
        $\size*{ \set*{j\in[\ell]_2\colon \mc{C}^{(j)}\neq\emptyset} }\leq\Gamma$; \label{item: conflict degree}  

		\item $\Delta_{j'}(\mc{C}^{(j)})\leq d^{j-j'-\eps}$ for all~$j\in[\ell]_2$ and~$j'\in[j-1]_2$;\label{item: conflict codegrees}
        
		\item $\size*{\set*{ f\in N^{(2)}_\mc{C}(e)\colon  v\in f }}\leq d^{1-\eps}$ for all~$e\in E(G)$ and $v\in V(G)$; \label{item: conflict j=2}
        
		\item $\size*{ N_\mc{C}^{(2)}(e)\cap N_\mc{C}^{(2)}(f)  }\leq d^{1-\eps}$ for all disjoint~$e,f\in E(G)$. \label{item: conflict double j=2}
	\end{enumerate}
\end{definition}

The goal of the conflict-free machinery is a $\mc{C}$-free matching with the additional properties which random edge sets that include the edges of~$G$ independently with probability $\size{M}/\size{E(G)}$ typically exhibit. In 2022, Glock, Joos, Kim, K\"
{u}hn, and Lichev~\cite{GJKKL22} showed that such a matching~$M$ can be chosen such that for a given sufficiently large edge set~$Z\subseteq E(G)$, we have~$\size{Z\cap M}\approx\size{Z} \cdot \size{M}/\size{E(G)}$, which is what would be expected if the edges of~$G$ were included in~$M$ independently at random with probability $\size{M}/\size{E(G)}$. To that end, we say that a $j$-uniform hypergraph~$\mc{Z}$ with vertex set~$E(G)$ whose edges are matchings of $G$ is a {\em test-system} and we call edges of~$\mc{Z}$ {\em tests}. However, we can only allow test systems that are well behaved in the sense that they can keep track of their properties during the iterative matching construction.

\begin{definition}[Trackable test systems]\label{def:trackable}
    Let $G$ be a graph. For integers~$j,d\geq 1$, a real~$\eps>0$ and a conflict system~$\mc{C}$ for~$G$ with $|C|\leq\ell$ for all $C\in\mc{C}$, we say that a~$j$-uniform test system~$\mc{Z}$ for~$G$ is~$(d,\eps,\mc{C})${\em-trackable} if the following holds.
	\begin{enumerate}[label=(Z\arabic*)]
		\item $\size*{\mc{Z}}\geq d^{j+\eps}$;\label{item: trackable size}
        
		\item $\Delta_{j'}(\mc{Z})\leq \size{\mc{Z}}/d^{j'+\eps}$ for all~$j'\in[j-1]$; \label{item: trackable degrees} %\TL{temp note: $\mc{Z}$ is $j$-uniform, so no need for $\mc{Z}^{(j)}$.}
        
		\item $\size{\mc{C}_e^{(j')}\cap \mc{C}_f^{(j')}}\leq d^{j'-\eps}$ for all~$e,f\in G$ with~$d_{\mc{Z}}(ef)\geq 1$ and all~$j'\in[\ell-1]$;\label{item: trackable neighborhood}
        
		\item $Z$ is~$\mc{C}$-free for all~$Z\in\mc{Z}$. \label{item: trackable no conflicts}
	\end{enumerate}
\end{definition}

We now recall the conflict-free matching machinery of Glock, Joos, Kim, K\"{u}hn, and Lichev~\cite{GJKKL22}, with test-systems. For clarity, we state only the special case needed in our proof.

\begin{thm}[Glock, Joos, Kim, K\"{u}hn, and Lichev~\cite{GJKKL22}]\label{thm:ConflictFree}
    For every sufficiently small~$\eps>0$, there exists~$d_0$ such that the following holds for all~$d\geq d_0$.
    Suppose~$r\geq 2$ is an integer and suppose~$\Gamma\geq 1$ and~$\mu\in(0,1/r]$ are reals such that~$1/\mu^{\Gamma r}\leq d^{\eps^2}$.
    Suppose~$G$ is a graph on~$n\leq \exp(d^{\eps^2/r})$ vertices with
    \((1-d^{-\eps})d\leq \delta(G)\leq\Delta(G)\leq d\). %, and \(\Delta_2(\mc{H})\leq d^{1-\eps}.\)
    Suppose that~$\mc{C}$ is a~$(d,r,\Gamma,\eps)$-bounded conflict system for~$G$, and that~$\mc{Z}$ is a $j$-uniform~$(d,\eps,\mc{C})$-trackable test system for~$G$ for some integer~$j\leq r$.
	Then, there exists a $\mc{C}$-free matching~$M\subseteq G$ of size~$(1-\mu)n/2$ with~$\size{\set{Z\in \mc{Z}\colon Z\subseteq M }}=(1\pm d^{-\eps/900})\parens{\size{M}/\size{G}}^j\size{\mc{Z}}$.
\end{thm}

We also require the following corollary of the polynomial concentration theorem of Kim-Vu~\cite{KV00}, included in a more general form in~\cite{DKPIII}, to prove concentration results on degrees and codegrees of various auxiliary hypergraphs used in the proof of~\cref{thm:maincoloring}.

\begin{cor}\label{cor:KimVu}
For all integers $k\ge 1$, there exists $n_0 \ge 1$ such that the following holds for all $p\in (0,1]$ and $n\ge n_0$. Let $H$ be a $k$-uniform multi-hypergraph on at most $n$ vertices, and $H_p$ be the subgraph of $H$ induced on the set $X$ where $X$ is obtained by choosing each vertex of $H$ independently with probability $p$. If $\kappa>0$ satisfies
\begin{equation}
    p^{i}\ge \frac{\Delta_i(H)}{\kappa}\cdot \log^{4k+2} n\tag{$*_i$}
\end{equation}
for all $i \in [k]$ and $e(H)\leq \kappa$, then with probability at least $1-e^{-\log^2 n}$,
\[e(H_p)\le p^k \cdot 2\kappa.\]
\end{cor}

Finally, we need the following deterministic lemma to increase a graph while preserving its maximum degree and ensuring that the resulting graph has sufficiently many edges. This will ensure that the test system constructed in the proof of~\cref{thm:maincoloring} is sufficiently large.

\begin{lem}\label{lem:supertestsystem}
Let $D\geq 2$ be an integer. If $H$ is a graph on $n\geq D+1$ vertices with $\Delta(H) \le D$, then there exists a supergraph $F$ of $H$ on $n$ vertices with $\Delta(F) \le D$ and $e(F)\geq Dn/4$.
\end{lem}

\begin{proof}
Let $F$ be a supergraph of $H$ on $n$ vertices with $\Delta(F) \leq D$ and, subject to that, $e(F)$ maximized. If every vertex of $F$ has degree at least $D/2$, then $e(F) \geq Dn/4$ as desired. Now assume there exists a vertex $v \in V(F)$ with degree strictly less than $D/2$. Hence, $d_{F}(v)\leq (D-1)/2 \leq (n-2)/2 = n/2 - 1$. If $v$ is non-adjacent to a vertex $w$ of degree strictly less than $D$, then $F'=F+vw$ contradicts the choice of $F$. Thus, every non-neighbor of $v$ has degree exactly $D$. Since $n\geq D$ and $d_{F}(v) < D/2$, we have that there are at least $n-1-d_{F}(v) \geq  n-1- (n/2-1) \geq n/2$ such vertices $w$. Hence the sum of degrees in $F$ is at least $(n/2)\cdot D$ and so $e(F)\geq Dn/4$ as desired. 
\end{proof}

\subsection{Proof of main technical result}\label{subsec:IndepPf}

Now we are prepared to prove~\cref{thm:maincoloring}.

\begin{lateproof}{thm:maincoloring}

    %%%%%%%%%%%%%%%%%%%%%%%%%%%
    %1. Parameters
    %%%%%%%%%%%%%%%%%%%%%%%%%%%
    Let $r\ge 3$ be an integer. Let $\beta,K > 0$ be real. Without loss of generality,  assume that $\beta\leq 1$. Let $\eps\in(0,1)$ be arbitrarily small, in particular such that~\cref{thm:ConflictFree} holds. Let $\mu:=\frac1r$ and let $\Lambda$ be large enough, in particular such that $\Lambda^2\cdot \eps^2\beta>80r^2\log r$ and $\Lambda\geq 32$. Let $d_0:=d_0(\eps)$ be given by~\cref{thm:ConflictFree}. Let $q$ be a sufficiently large integer,  and 
    \[\chi:=\ceil*{\Lambda q^2},
    \qquad d:=\floor*{q^{\beta/10}}>d_0,
    \qquad \Gamma:=\frac{8r}{\Lambda^2}\cdot\log q.\]
    Let $D$ be large enough compared to $q$. Let $H$ be an $r$-bounded hypergraph on $N$ vertices with $N\leq KqD$, whose edges have size at least two. Observe that we can assume that $N\geq 16D$, as otherwise the statement is trivially true by taking $c<1/16$. Assume that for all integers $2\leq \ell<j\leq r$ we have
    \[\begin{cases}
        \Delta_{1}(H^{(2)})&\leq D\\
        \Delta_{1}(H^{(j)})&\leq D^{j-1}\cdot\log q\\
        \Delta_{\ell}(H^{(j)})&\leq D^{j-\ell}\cdot q^{-\beta}.
    \end{cases}\]
    Recall that we aim to prove that $\alpha(H)\geq c\cdot\frac{N}{D}$ for some constant $c=c(r,\beta)$.\bigskip

    %%%%%%%%%%%%%%%%%%%%%%%%%%%
    %2. Random subgraph H'
    %%%%%%%%%%%%%%%%%%%%%%%%%%%  
    Let $p:=q^2/D$, and let $H'\subseteq H$ be the induced subgraph of $H$ where every vertex is selected with probability $p$, independently of all other vertices. Standard concentration inequalities show that with high probability (for $D$ large enough, as $N\geq 16D$), we have
    \begin{equation}\label{eq:Bounds_VHPrime}
       8q^2\leq \frac{q^2 N}{2D}\leq m:=v(H')\leq \frac{2q^2N}{D}\leq 2Kq^3,
    \end{equation}
    and
    \begin{equation}\label{eq:Bounds_DeltaHPrime}
        \Delta_{1}(H'^{(2)})\leq 2q^2,\qquad
        \Delta_{1}(H'^{(j)})\leq 2q^{2(j-1)}\cdot\log q,\qquad
        \Delta_{\ell}(H'^{(j)})\leq 2q^{2(j-\ell)-\beta}.
    \end{equation}
    While the bounds on $v(H')$ and $\Delta_1(H'^{(2)})$ follow from standard application of the Chernoff and union bounds, the upper bounds on $\Delta_{1}(H'^{(j)})$ and $\Delta_{\ell}(H'^{(j)})$ are direct application of tail estimates for randomly induced subhypergraphs, applied to the links of $H$ (see e.g. Janson and Ruci\'{n}ski~\cite[Section 2]{JR2011}). We fix one such a realization $H'$.
    
    %%%%%%%%%%%%%%%%%%%%%%%%%%%
    %3. near-twins
    %%%%%%%%%%%%%%%%%%%%%%%%%%%  
    For $j \in [r]$, say vertices $u,v\in V(H')$ are \emph{$j$-near twins} if there are at least $q^{2(j-1)-\beta/2}$ distinct $(j-1)$-sets $T$ such that $T\cup\{u\} \in E(H)$ and $T \cup \{v\} \in E(H)$, that is, if 
    \[\size{{H'}_u^{(j-1)}\cap {H'}_v^{(j-1)}}\geq q^{2j-2-\beta/2}.\]
    Let $W_j$ be the graph of $j$-near twins, that is, the graph with $V(W_j):=V(H')$ and $uv\in E(W_j)$ if and only $u,v$ are $j$-near twins. Let $W:=\bigcup_{j=3}^rW_j$. We first derive an upper bound on $\Delta(W)$. For every $u\in V(H')$ and every $j\in\{3,\ldots,r\}$, observe that, for every $v\in N_{W_j}(u)$ we have $\size{{H'}_u^{(j-1)}\cap {H'}_v^{(j-1)}}\geq q^{2j-2-\beta/2}$. Therefore
    \[d_{W_j}(u)\leq \frac{1}{q^{2j-2-\beta/2}}\cdot\sum_{v\in N_{W_j}(u)} \size{{H'}_u^{(j-1)}\cap {H'}_v^{(j-1)}}.\]
    Meanwhile, by~\cref{eq:Bounds_DeltaHPrime},
    \[\sum_{v\neq u}\size{{H'}_u^{(j-1)}\cap {H'}_v^{(j-1)}}
    =\sum_{T\in {H'}_u^{(j-1)}} \parens*{\deg_{{H'}^{(j)}}(T)-1}\leq \Delta_1({H'}^{(j)})\cdot \Delta_{j-1}({H'}^{(j)})
    \leq 4q^{2j-\beta}\cdot\log q,
    \]
    and therefore
    \begin{equation}\label{eq:Bound_DeltaW}
        d_{W_j}(u)\leq4q^{2-\beta/2}\cdot \log q, \quad \text{ and }\qquad
        \Delta(W)\leq4rq^{2-\beta/2}\cdot \log q\ll q^2.    
    \end{equation}
    
    %%%%%%%%%%%%%%%%%%%%%%%%%%%%%
    %4. Supergraph F and F':=F\W
    %%%%%%%%%%%%%%%%%%%%%%%%%%%%%
    We apply~\cref{lem:supertestsystem} to ${H'}^{(2)}$ with parameter $4q^2$, as ${H'}^{(2)}$  is a graph on $m\geq 8q^2$ vertices, with maximum degree at most $2q^2$. We obtain a supergraph $F\supseteq {H'}^{(2)}$ on $m$ vertices with $\Delta(F)\leq 4q^2$, at least $mq^2$ edges, and at most $2mq^2$.
    
    Let $F':=F\setminus W$. Since $\Delta(W)\leq4rq^{2-\beta/2}\cdot \log q$, for $q$ large enough we obtain that
    \begin{equation}\label{eq:Bound_DeltaF}
    \frac{mq^2}{2}\leq e(F')\leq 2mq^2,\qquad
    \Delta(F')\leq 4q^2,\qquad
    {H'}^{(2)}\setminus W\subseteq F'\setminus W,\qquad
    F'\cap W=\emptyset
    \end{equation}

    %%%%%%%%%%%%%%%%%%%%%%%%%%%%%%%%%%%%%%%
    %5. Sparse color palette and Graph G, 
    %%%%%%%%%%%%%%%%%%%%%%%%%%%%%%%%%%%%%%% 
    We now encode $d$ partial colorings of $H'$ as a matching in an auxiliary graph $G$. The vertex $(v,i,0)$ in $V(G)$ will represent the vertex $v$ in the $i$th coloring, while $(v,c,1)$ in $V(G)$ will represent the use of the color $c\in L(v)$. Therefore, selecting the edge $e^{(i)}_{v,c}$ in a matching of $G$ will mean assigning color $c$ to $v$ in the $i$th coloring.   
    
    For each $v\in V(H')$ let $L(v)$ be a (uniformly) random $d$-subset of $[\chi]$. For each $v\in V(H')$, let $G_v\cong K_{d,d}$ with vertex set $(A_v,B_v)$ where
    \[A_v:= \set*{(v,i,0)\colon i\in [d]},\qquad B_v:= \set*{ (v,c,1)\colon c\in L(v)}.\]
    Let $G$ be the disjoint union of copies of $G_v$, that is, $G:=\bigsqcup_{v\in V(H')}G_{v}$, and observe that we have 
    \[|V(G)|=2d\cdot |V(H')|=2dm,\quad |E(G)|=m d^2.\]

    For every $v\in V(H')$, $i\in[d]$, and $c\in L(v)$, let $e^{(i)}_{v,c}$ be the edge $(v,i,0)(v,c,1)$ of $G$. Let 
    \[E_1:=\set*{  \big\{e^{(i)}_{v_1,c},\ldots,e^{(i)}_{v_j,c}\big\}
    \colon \{v_1,\ldots,v_j\}\in E(H'),\ j\geq 3,\ c\in \bigcap_{\ell=1}^jL(v_\ell),\ i\in [d]},
        \]
    and 
    \[E_2: = \set*{\big\{e^{(i)}_{u,c}, e^{(i)}_{v,c}\big\}
    \colon \{u, v\} \in W, c \in L(u)\cap L(v), i \in [d]}.\]
    Let $\mc{C}$ be the conflict system for $G$ defined as $E(\mc{C}) = E_1 \cup E_2$. Observe that the conflicts in $E_1$ forbid monochromatic edges of $H'$ of size at least three, while those in $E_2$ forbid monochromatic pairs of near twins.   
    For an edge $e\in E(G)$ and an integer $\ell\in [r]$, we write
    \[\mc{C}_e^{\ell}:=\set*{C\setminus e\colon~e\in C\in\mc{C},|C|=\ell+1}.\]

    In order to apply the conflict-free~\cref{thm:ConflictFree}, we must first verify the degree and codegree bounds of~\cref{def:BoundedConflict}.~\Cref{clm:OneDegree} follows from Chernoff's inequality and a union bound, while~\cref{clm:LargeDegree,clm:LargeCoDegree}
    are applications of the Kim-Vu~\cref{cor:KimVu}. To keep the main argument focused, we defer their proofs to~\cref{app:concentration}.

    \begin{claim}\label{clm:OneDegree}
        For the near twins conflict, we have
        \[\mathbb{P}\brackets*{\Delta_{1}(\mc{C}^{(2)})\geq d^{1-\varepsilon}}\leq 1/10.\]
    \end{claim}

    \begin{claim}\label{clm:LargeDegree}
        We have
        \[\mathbb{P}\brackets*{\exists~j\in\{3,\ldots,r\}\colon~\Delta_{1}(\mc{C}^{(j)})\geq \frac{4}{\Lambda^{j-1}}\cdot d^{j-1}\log q }\leq\frac{1}{10}.\]
    \end{claim}
    
    \begin{claim}\label{clm:LargeCoDegree}
        We have,
        \[\mathbb{P}\brackets*{\exists~2\leq j'<j\leq r\colon~\Delta_{j'}(\mc{C}^{(j)})\geq d^{j-j'-\varepsilon}}\leq\frac{1}{10}.\]
    \end{claim}

    Using~\crefrange{clm:OneDegree}{clm:LargeCoDegree}, we now prove that the conflict-system is bounded.
    \begin{claim}\label{clm:ConflictBounded}
        The conflict-system $\mc{C}$ is $(d,r,\Gamma,\eps)$-bounded with probability at least $3/5$.
    \end{claim}
    \begin{proofclaim}
        Property C1) is easily seen to be true because $H$ is $r$-bounded. For Property C2), observe that, for $q$ large enough compared to $r$,~\cref{clm:OneDegree,clm:LargeDegree} imply that, with probability at least~$4/5$ we have 
        \[\sum_{j\in[r]} \frac{\Delta_{1}(\mc{C}^{(j)})}{d^{j-1}} \leq 
        d^{-\varepsilon}+ \sum_{j=3}^{r} \frac{4}{\Lambda^{j-1}}\cdot\log q 
        \leq 1 + \frac{4r}{\Lambda^2}\cdot\log q
        \leq \frac{8r}{\Lambda^2}\cdot\log q = \Gamma.\]
        We furthermore have
        \(\size*{ \set*{j\in[r]_2\colon \mc{C}^{(j)}\neq\emptyset} }\leq r-1\leq \Gamma, \) hence Property C2) is satisfied with probability at least~$4/5$. Property C3) is a direct consequence of~\cref{clm:LargeCoDegree}, and therefore holds with probability at least $9/10$.

        To prove Property C4), let $e=e_{u,c}^{(i)}$ be an edge of $G$ between vertices $(u,i,0)$ and $(u,c,1)$. Every $2$-conflict neighbor of $e$ has the form $f=e_{v,c}^{(i)}$ with $uv\in W$. Fix a vertex $z\in V(G)$. If $z=(w,i',0)$, then $z\in f$ forces $v=w$ and $i=i'$, hence $f$ is uniquely determined by $e$ and $z$; If $z=(w,c',1)$, then again $z\in f$ forces $v=w$ and $c=c'$, hence $f$ is again uniquely determined by $e$ and $z$. It follows that, deterministically we have 
        \[\size*{\set*{ f\in N^{(2)}_\mc{C}(e)\colon  v\in f }}\leq 1\leq d^{1-\eps},\] as desired. Finally, Property C5) holds with probability at least $9/10$,as it follows from~\cref{clm:OneDegree} that, for any two distinct edges $e,f\in E(G)$,
        \[\size*{ N_\mc{C}^{(2)}(e)\cap N_\mc{C}^{(2)}(f)  }
        \leq\size*{ N_\mc{C}^{(2)}(e) }
        \leq \Delta_{1}(\mc{C}^{(2)})
        \leq d^{1-\varepsilon}.
        \]
        Combining these properties, we obtain that $\mc{C}$ is $(d,r,\Gamma,\eps)$-bounded with probability at least $3/5$, as desired.
    \end{proofclaim}

    The conflict system $\mc{C}$ defined above prevents monochromatic edges of size at least three and monochromatic pairs of near twins. We cannot directly forbid the remaining monochromatic $2$-edges with the current (co)degree requirements of the conflict-free machinery. Instead, we record them through a test system and later delete one endpoint from every test. This alteration will remove every monochromatic $2$-edge not already forbidden by a conflict. Let $\mc{Z}$ be the test system for $G$ defined by
    \[E(\mc{Z}):=\set*{ 
    \big\{e^{(i)}_{u,c},e^{(i)}_{v,c}\big\}
    \colon u \neq v,\ \{u, v\} \in F',\ c \in L(u)\cap L(v),\ i\in[d]}.\]  
    
    We first estimate the size of the test system $\mc{Z}$. The lower bound will be used to verify condition (Z1) of~\cref{def:trackable}, while the upper bound will control the number of vertices deleted in the final alteration. A second-moment calculation followed by Chebyshev's inequality gives the following estimate; the details are again deferred to~\cref{app:concentration}.
    \begin{claim}\label{clm:SizeConflict}
        We have
        \[\mathbb{P}\brackets*{\frac{m\cdot  d^3}{8\Lambda}\leq |\mc{Z}|\leq \frac{4m\cdot  d^3}{\Lambda}}\geq \frac{9}{10}.\]
    \end{claim}    
    
    We now combine various estimates above to verify that $\mc{Z}$ is trackable.    
    \begin{claim}\label{clm:TestTrackable}
        The $2$-uniform test-system $\mc{Z}$ is $(d,\eps,\mc{C})$-trackable with probability at least $7/10$.
    \end{claim}    
    \begin{proofclaim}
        By~\cref{clm:SizeConflict}, with probability at least $9/10$ we have \[|\mc{Z}|\geq \frac{m\cdot  d^3}{8\Lambda}.\]
        Using $m\geq 8q^2$, for $q$ large enough, we have
        \[\frac{m\cdot  d^3}{8\Lambda} \geq \frac{q^2\cdot  d^3}{\Lambda}\geq d^{2+\varepsilon},\qquad\text{and}\qquad
        \frac{m\cdot d^3}{8\Lambda} \geq d^{1+\varepsilon}\cdot \frac{q^2\cdot  d^{2-\varepsilon}}{\Lambda}\geq d^{1+\varepsilon}\cdot 4q^2,\]
        while~\cref{eq:Bound_DeltaF} implies
        \[\Delta_1(\mc{Z})\leq\Delta(F')\leq 4q^2.\]
        % \[\frac{|\mc{Z}|}{d^{1+\varepsilon}}\geq \Delta_1(\mc{Z}),\]
        Therefore it follows that Properties Z1) and Z2) hold with probability at least $9/10$. Observe that, by~\cref{eq:Bound_DeltaF} we have $F\cap W=\emptyset$, hence we trivially have that $Z$ is~$\mc{C}$-free for all~$Z\in\mc{Z}$. It remains to prove that Property Z3) holds with sufficient large probability. We do so via a standard application of~\cref{cor:KimVu}. 

        Recall that we aim to show that $\size{\mc{C}_e^{(j')}\cap \mc{C}_f^{(j')}}\leq d^{j'-\eps}$ for all~$e,f\in G$ with~$d_{\mc{Z}}(ef)\geq 1$ and all~$j'\in[r-1]$;

        If $j'=1$, then, similarly to the proof of Property C5) in~\cref{clm:ConflictBounded}, it follows from~\cref{clm:OneDegree} that,  with probability at least $9/10$, for every $e,f\in E(G)$
        \[ \size{\mc{C}_e^{(j')}\cap \mc{C}_f^{(j')}} 
        = \size*{ N_\mc{C}^{(2)}(e)\cap N_\mc{C}^{(2)}(f)  } 
        \leq\size*{ N_\mc{C}^{(2)}(e) }
        \leq \Delta_{1}(\mc{C}^{(2)})        
        \leq d^{1-\varepsilon}.\]
        
        Assume now that $2\leq j'\leq r-1$. Fix $e,f\in E(G)$ with $d_{\mc{Z}}(ef)\geq1$. Since $\mc{Z}$ is $2$-uniform, $ef$ is an edge of $\mc{Z}$ and by definition of the test-system, there exist distinct vertices $u,v\in V(H')$, a color $c\in[\chi]$, and an integer $i\in[d]$ such that
        \[ e=e_{u,c}^{(i)}, \quad f=e_{v,c}^{(i)}, \quad uv\in E(F').\]

        We define the $j'$-uniform hypergraph $J=J_{u,v,j'}$ by 
        \[J:={H'}_u^{(j')}\cap {H'}_v^{(j')}.\]
        Let $\rho:=d/\chi$ and, for each color $c\in[\chi]$, let \[X_c:=\set*{w\in V(H')\colon c\in L(w)}.\] Note that for every $c \in [\chi]$, every vertex of $V(H')$ is in $X_c$ independently with probability $\rho$. Now let $J_{c, \rho}$ be the random subhypergraph of $J$ induced on the set 
        \[X_c:=\set*{w\in V(H')\colon c\in L(w)}.\]
        For each color $c\in[\chi]$, every vertex of $V(H')$ is in $X_c$ independently with probability $\rho$. We claim that $e(J_{c, \rho})$ represents the quantity we are trying to estimate. Indeed, for every $S\in E(J_{c, \rho})$, let $Q_S:=\set*{e^{(i)}_{w,c}:w\in S}$. Since $S\cup\{u\},S\cup\{v\}\in E(H')$, we have that both $Q_S\cup\{e\}$ and $Q_S\cup\{f\}$ are conflicts, and therefore 
        \(Q_S\in\mathcal C_e^{(j')}\cap\mathcal C_f^{(j')}.\)
        Conversely, every member of this intersection is uniquely of the form $Q_S$ for some $S\in E(J_{c, \rho})$. Therefore
        \[\size{\mc{C}_e^{(j')}\cap \mc{C}_f^{(j')}}  = e(J_{c, \rho}).\]
        Let 
        \[\kappa:=2q^{2j'-\beta/20}.\]
        As $u,v$ are such that $uv\in E(F')$, we know that $u,v$ are not $(j'+1)$-near twins, hence
        \[e(J)\leq q^{2j'-\beta/2}\leq \kappa.\]
        We aim to apply~\cref{cor:KimVu} to the $j'$-uniform hypergraph $J$. We check first that $(*_\ell)$ holds for every $\ell\in[j'-1]$. Indeed, for any such $\ell$, by~\cref{eq:Bounds_DeltaHPrime} we have
        \[\Delta_{\ell}(J)\leq\Delta_{\ell+1}({H'}^{(j'+1)})\leq 2q^{2(j'-\ell)-\beta},\]
        hence, with $d=\floor{q^{\beta/10}}$ and $\chi=\ceil{\Lambda q^2}$, for $q$ large enough
        \[\frac{\Delta_{\ell}(J)}{\kappa}\cdot \log^{4j'+2}m\leq
        q^{-2\ell-\beta+\beta/20} \cdot \log^{4j'+2}m\leq\parens*{\frac{d}{\chi}}^\ell = \rho^{\ell}.\]
        If $\ell=j'$, we have $\Delta_{\ell}(J)\leq 1$, and for $q$ large enough,
        \[\frac{\Delta_{\ell}(J)}{\kappa}\cdot \log^{4j'+2}m
        \leq \frac12 \cdot q^{-2j'+\beta/20} \cdot \log^{4j'+2}(q^4)
        \leq q^{-2j'}\cdot q^{\beta/10} 
        \leq\parens*{\frac{d}{\chi}}^\ell = \rho^{\ell},\]
        hence $*_{j'}$ also holds. Therefore, by~\cref{cor:KimVu}, we have 
        \[e(J_{c, \rho})\geq 2\rho^{j'}\kappa\]
        with probability at most $e^{-\log^2m}$. Observe finally that
        \[2\rho^{j'}\kappa\leq \frac{4}{\Lambda^{j'}}\cdot \frac{d^{j'}}{q^{2j'}}\cdot q^{2j'-\beta/20}\leq d^{j-\varepsilon},\]
        as $\varepsilon$ is small enough compared to $\beta$. We conclude that, for $q$ large enough
        \[\mathbb{P}\brackets*{\exists j\in \{2,\ldots,r-1\},~e,f\in G\text{ with }d_{\mc{Z}}(ef)\geq 1\colon \size{\mc{C}_e^{(j')}\cap \mc{C}_f^{(j')}} \geq d^{j'-\varepsilon}} \leq r \cdot  (md^2)^2\cdot e^{-\log^2 m}\leq \frac{1}{10}.\]
        Therefore Property (Z3) holds with probability at least $8/10$, hence $\mc{Z}$ is $(d,\eps,\mc{C})$-trackable with probability at least $7/10$, as desired.             
    \end{proofclaim}
        
    We obtain from~\cref{clm:ConflictBounded,clm:TestTrackable}, that, with probability at least $3/10$, $\mc{C}$ is $(d,r,\Gamma,\eps)$-bounded and $\mc{Z}$ is $(d,\eps,\mc{C})$-trackable. Fix such an outcome. By definition $G$ is d-regular, with $|V(G)|$ being polynomial in $d$, therefore, for $q$ large enough we have $|V(G)|\leq \exp(d^{\varepsilon^2/r})$. To apply~\cref{thm:ConflictFree}, it remains to show that $1/ \mu^{\Gamma r}\leq d^{\eps^2}$. Recall that $\mu=\frac1r$, $\Gamma=\frac{8r}{\Lambda^2}\cdot\log q$, and $d=\floor*{q^{\beta/10}}$, therefore
    \[d^{\eps^2}\cdot \mu^{\Gamma r} \geq q^{\eps^2\beta/10}\cdot \mu^{8r^2\log q/\Lambda^2},\]
    and
    \[\log_q\parens*{d^{\eps^2}\cdot \mu^{\Gamma r}} 
    \geq \frac{\eps^2\beta}{10}+ \frac{8r^2}{\Lambda^2}\log \mu
    =    \frac{\eps^2 \cdot \beta}{10}- \frac{8r^2 \cdot \log r}{\Lambda^2}>0,\]
    as $\Lambda$ is such that $\Lambda^2\cdot \eps^2 \cdot \beta>80r^2\cdot \log r$. It follows that~$1\leq d^{\eps^2}\cdot \mu^{\Gamma r}$, as desired. We then apply the conflict-free theorem,~\cref{thm:ConflictFree}, and we obtain a $\mc{C}$-free matching~$M\subseteq G$ with
    \[\size{M} = \frac{1-\mu}{2}\cdot|V(G)|,\qquad \size{\set{Z\in \mc{Z}\colon Z\subseteq M }}=(1\pm d^{-\eps/900})\cdot\frac{\size{M}^2}{\size{E(G)}^2}\cdot\size{\mc{Z}}.\]

    Recall that $v(G)=2dm$ and that $e(G)=md^2$, hence
    \[\size{M} = (1-\mu)dm,\qquad \size{\set{Z\in \mc{Z}\colon Z\subseteq M }}=(1\pm d^{-\eps/900})\cdot\frac{(1-\mu)^2}{d^2}\cdot\size{\mc{Z}}.\]
    With $|\mc{Z}|\leq \frac{4m \cdot  d^3}{\Lambda}$, we have 
    \[\size{\set{Z\in \mc{Z}\colon Z\subseteq M }}
    \leq (1\pm d^{-\eps/900})\cdot(1-\mu)^2\cdot\frac{4m\cdot d}{\Lambda}
    \leq \frac{8m\cdot d}{\Lambda}.\]

    Recall that we view the matching $M$ as a family of $d$ partial colorings of $V(H')$, where the edge $e^{(i)}_{v,c}\in M$ represent the assignment of the color $c$ to $v$ in the $i$th coloring. For each $i\in[d]$, let $a_i$ be the number of vertices colored by $M$ in the $i$th-coloring, that is, edges of the form $e_{v,c}^{(i)}$ in $M$. Further let $b_i$ be the number of test-edges in $M$ involving the $i$th-coloring. By definition we have
    \[\sum_{i=1}^d a_i=|M|=(1-\mu)md,\qquad\sum_{i=1}^db_i\leq \frac{8m\cdot d}{\Lambda},\]
    hence, with $\Lambda\geq 32$ and $\mu=1/r\leq 1/2$,
    \[\sum_{i=1}^d (a_i-b_i)\geq (1-\mu)\cdot m\cdot d - \frac{8m\cdot d}{\Lambda} \geq \frac{1-\mu}{2}\cdot m\cdot d.\]
    Therefore there exists $i\in [d]$ such that
    \[(a_i-b_i)\geq \frac{1-\mu}{2}\cdot m.\]
    In that $i$th coloring, we delete one end-point of any monochromatic edge in $F'$, hence deleting at most $b_i$ vertices, hence we still have at least $\frac{1-\mu}{2}m$ colored vertices. We claim that this is a proper partial coloring of $H'$, that is every color class is independent; indeed:
    \begin{enumerate}
        \item Monochromatic $j$-edges with $j\geq 3$ are forbidden by the conflicts.
        \item Monochromatic $2$-edges that are near-twins are also in conflict, hence forbidden.
        \item Monochromatic $2$-edges that are not near-twins are in $F'$, hence deleted by alteration.
    \end{enumerate}
    At least one color class has size at least $\frac{(1-\mu)m}{2\chi}$. Recall that $m\geq \frac{q^2N}{2D}$ and $\chi\leq 2\Lambda q^2$, hence
    \[\alpha(H)\geq\alpha(H')\geq \frac{1-\mu}{8\Lambda}\cdot \frac{N}{D},\]
    as desired by letting $c(r,\beta):=\frac{1-\mu}{8\Lambda}$.
\end{lateproof}

\vspace{1em}

\noindent \textbf{Declaration of absence of AI use:} This project was initiated in early 2025 and the primary mathematical ideas were developed by the fall of 2025. No generative AI tools were used in the research or preparation of this article.

\bibliographystyle{plain}
\bibliography{bibliography}

%%%%%%%%%%%%%%%%%%%%%%%%%%%%%%%%%%%%%%%
%%%%% Appendix 
%%%%%%%%%%%%%%%%%%%%%%%%%%%%%%%%%%%%%%%

\appendix

\section{Concentration inequalities}\label{app:concentration}

In this appendix, we prove the concentration estimates from the proof of~\cref{thm:maincoloring}.~\Cref{clm:OneDegree} follows from Chernoff's inequality and a union bound;~\cref{clm:LargeDegree,clm:LargeCoDegree} follow from~\cref{cor:KimVu}; and~\cref{clm:SizeConflict} follows from Chebyshev's inequality.

\subsection{Maximum 1-degree}

\begin{proof}[Proof of~\cref{clm:OneDegree}]
    For a vertex $u\in V(H')$ and a color $c\in[\chi]$, let 
    \[X_{u,c}=\size*{\set*{v\in N_W(u)\colon c\in L(v)}},\]
    hence the number of near-twins of $u$ with the color $c$ available in the sparsified list. Every $2$-conflict containing the edge $e^{(i)}_{u,c}$ is of the form $\{e^{(i)}_{u,c},e^{(i)}_{v,c}\}$ with $v\in N_W(u)$ and $c\in L(v)$. Therefore, if we can prove that $X_{u,c}\leq d^{1-\varepsilon}$ for all choices of $u,c$ then we are done.
    By the independent random choice of the lists $L(v)$, we have $X_{u,c}\sim{\rm Bin}(d_W(u),\frac{d}{\chi})$, hence
    \[ \mathbb{E}X_{u,c} \leq \Delta(W)\frac{d}{\chi} 
    \leq 4rq^{2-\beta/2}\cdot\log q \cdot\frac{d}{\Lambda q^2}
    \leq \frac{4r}{\Lambda} d^{-4}\cdot \log q,
    \]
    where we used $d\leq q^{\beta/10}$. Let $B,\delta>0$ such that $B:=(1+\delta) \mathbb{E}X_{u,c} = d^{1-\varepsilon}$. By the Chernoff bound, 
    \[\mathbb{P}\brackets*{X_{u,c}\geq d^{1-\varepsilon}}\leq \parens*{\frac{e^\delta}{(1+\delta)^{1+\delta}}}^{\mathbb{E}X_{u,c}}\leq \parens*{\frac{e\cdot \mathbb{E}X_{u,c}}{B}}^B\leq 
    e^{- d^{1-\varepsilon}}\]
    By union bound over the at most $m\chi=O(q^5)$ pairs $(u,c)$, we obtain that, for $q$ large enough
    \[\mathbb{P}\brackets*{\Delta_{1}(\mc{C}^{(2)})\geq d^{1-\varepsilon}}\leq e^{- d^{1-\varepsilon}}\cdot m\chi \leq \frac{1}{10},\]
    as desired.
\end{proof}

\subsection{Maximum j-degrees}

\begin{proof}[Proof of~\cref{clm:LargeDegree}]
    Fix $j\in\{3,\ldots,r\}$, a vertex $v\in V(H')$, and a color $c\in[\chi]$. Let $J=J(v,j)$ be the $(j-1)$-uniform hypergraph defined by $V(J)=V(H')\setminus\{v\}$ and
        \[E(J)={H'}^{(j-1)}_v =\set*{ S\in\binom{V(H')\setminus\{v\}}{j-1}\colon S\cup\{v\}\in {H'}^{(j)}}.\]
    Let $\rho=d/\chi$ and let $J_{\rho}$ be the random subhypergraph of $J$ induced on the set
    \[X_c :=\set*{u\in V(H')\colon c\in L(u)}.\]
    For each color $c\in[\chi]$, every vertex of $V(H')$ is in $X_c$ independently with probability $\rho$. If $c\in L(v)$, then, for every $i\in [d]$ we have
    \[d_{\mc{C}^{(j)}}(e_{v,c}^{(i)})=e(J_{c, \rho}).\]
    Let 
    \[\kappa:=2q^{2(j-1)}\cdot\log q.\]
    We have
    \[e(J)=d_{{H'}^{(j)}}(v)\leq \Delta_{1}({H'}^{(j)})\leq 2q^{2(j-1)}\cdot\log q=\kappa.\]
    We apply~\cref{cor:KimVu} to $J$ (which is $(j-1)$-uniform). For every $\ell\in[j-2]$
    \[\Delta_\ell(J)\leq \Delta_{\ell+1}({H'}^{(j)})\leq 2q^{2(j-\ell-1)}\cdot q^{-\beta},\]
    We therefore can check that $(*_\ell)$ holds, for $q$ large enough       
    \begin{align*}
    \frac{\Delta_\ell(J)}{\kappa}\cdot \log^{4(j-1)+2}m
    \leq q^{-2\ell -\beta}\cdot \frac{ \log^{4(j-1)+2}m}{\log q}
    \leq q^{-2\ell -\beta}\cdot 4^{^{4(j-1)+2}}\log^{4(j-1)+1}q
    \leq \parens*{\frac{d}{\chi}}^\ell
    =\rho^\ell
    \end{align*}
    where we use $m\leq q^4$, $d=\floor*{q^{\beta/10}}$ and $\chi=\ceil*{\Lambda q^2}$.
    For $\ell=j-1$, as $\Delta_{j-1}(J)\leq 1$, we obtain
    \[\frac{\Delta_{j-1}(J)}{\kappa}\cdot \log^{4(j-1)+2}m
    \leq \frac{\log^{4j-2}(q^4)}{2q^{2(j-1)}\log q}
    \leq\parens*{\frac{d}{\chi}}^{j-1}
    =\rho^{j-1}.\]
    Observe finally that
    \[2\rho^{j-1}\kappa = 4q^{2(j-1)}\cdot\log q \cdot\parens*{\frac{d}{\chi}}^{j-1}\leq \frac{4}{\Lambda^{j-1}}\cdot d^{j-1}\log q.\]
    Therefore, by~\cref{cor:KimVu}, we have
    \[e(J_{\rho})\geq \frac{4}{\Lambda^{j-1}}\cdot d^{j-1}\log q 
    \geq 2\rho^{j-1}\kappa\]
    with probability at most $e^{-\log^2m}$. It follows that
    \begin{align*}
        % \mc{P}:=~&
        \mathbb{P}\brackets*{\Delta_{1}(\mc{C}^{(j)})\geq \frac{4}{\Lambda^{j-1}}\cdot d^{j-1}\log q }
        \leq~&\mathbb{P}\brackets*{\exists 3\leq j\leq r,~v\in V(H'),~c\in[\chi],~i\in[d]\colon e(J(v,j)_{\rho})\geq \frac{4}{\Lambda^{j-1}}\cdot d^{j-1}\log q}\\
        \leq~&rm\chi d\cdot e^{-\log^2m}\\
        \leq~&O(q^6)\cdot e^{-\log^2q}\\
        \leq~&\frac{1}{10},
    \end{align*}
    as desired.
\end{proof}

\subsection{Maximum codegrees}
\begin{proof}[Proof of~\cref{clm:LargeCoDegree}]
    Fix integers $2\leq j'<j\leq r$, and let $k:=j-j'$. Fix some distinct vertices $v_1,\ldots,v_{j'}\in V(H')$, a color $c\in[\chi]$, and $i\in[d]$.  Let $T:=\{v_1,\ldots,v_{j'}\}$ and let $J=J(T,j)$ be the $k$-uniform hypergraph defined by $V(J):=V(H')\setminus T$ and 
    \[ E(J):= \set*{ S\in\binom{V(H')\setminus T}{k}\colon T\cup S\in H'^{(j)} }.\]
    Let $\rho=d/\chi$ and let $J_{\rho}$ be the random subhypergraph of $J$ induced on the set $X_c :=\set*{u\in V(H')\colon c\in L(u)}$. If $c\in L(v)$, then, for every $i\in [d]$ we have
    \[d_{\mc{C}^{(j)}}\parens*{\set*{e_{v,c}^{(i)}\colon v\in T}}=e(J_{c, \rho}).\]
    Let 
    \[\kappa:= 2q^{2k-\beta/20}\]
    We have
    \[e(J)=d_{{H'}^{(j)}}(T)\leq 2q^{2k-\beta} \leq \kappa.\]   
    We apply~\cref{cor:KimVu} to the $k$-uniform hypergraph $J$ . For every $\ell\in[k-1]$,
    \[ \Delta_\ell(J) \leq \Delta_{j'+\ell}(H'^{(j)}) \leq 2q^{2(k-\ell)-\beta}.\]        
    Using $\kappa\geq 2q^{2k-\beta}$, we now check that $(*_\ell)$ holds as we have 
    \[ \frac{\Delta_\ell(J)}{\kappa} \log^{4k+2}m \leq q^{-2\ell}\log^{4k+2}m
    \leq \parens*{\frac d\chi}^\ell
    =\rho^\ell,
    \]
    where we use $m\leq q^4$, $d=\floor*{q^{\beta/10}}$ and $\chi=\ceil*{\Lambda q^2}$. As $J$ is a simple $k$-uniform hypergraph, we have $\Delta_k(J)\leq 1$ and therefore, for $q$ large enough we have
    \[ \frac{\Delta_k(J)}{\kappa} \log^{4k+2}m \leq q^{-2k+\beta/20}\log^{4k+2}m
    \leq \parens*{\frac d\chi}^k
    =\rho^k,
    \]        
    where we used $\beta/20\ll \beta k/10$, hence $(*_k)$ also holds. Therefore, by~\cref{cor:KimVu}, we have
    \[e(J_{\rho})\geq 2\rho^{k}\kappa.\]
    with probability at most $e^{-\log^2m}$. Observe finally that
    \[2\rho^{k}\kappa \leq 2\cdot\parens*{\frac{q^{\beta/10}}{\Lambda q^2}}^k\cdot 2q^{2k-\beta/20}\leq\frac{4}{\Lambda^k}\cdot q^{\beta/10(k-1/2)}\leq d^{k-1/2}\leq d^{j-j'-\varepsilon},\]
    as $\varepsilon>0$ is arbitrarily small. It follows that
    \begin{align*}
        % \mc{P}:=~&
        \mathbb{P}\brackets*{\Delta_{j'}(\mc{C}^{(j)})\geq d^{j-j'-\varepsilon}}
        \leq~&\mathbb{P}\brackets*{\exists 2\leq j'<j\leq r,~T\in \binom{V(H')}{j'},~c\in[\chi],~i\in[d]\colon e(J(T,j)_{\rho})\geq 2\rho^{k}\kappa}\\
        \leq~&r^2m^{r}\chi d\cdot e^{-\log^2m}\\
        \leq~&O(q^{3(r+1)})\cdot e^{-\log^2q}\\
        \leq~&\frac{1}{10},
    \end{align*}
    as desired.    
\end{proof}

\subsection{Size of the test system}
\begin{proof}[Proof of~\cref{clm:SizeConflict}]
    We aim to prove upper and lower bounds on $|\mc{Z}|$, the size of our test-system. Observe that, 
    \[|\mc{Z}|=d\cdot\sum_{uv\in F'}|L(u)\cap L(v)|,\]
    hence
    \begin{align*}
    \mathbb{E}|\mc{Z}|
    =d\cdot \sum_{uv\in F'} \mathbb{E}[\size{L(u)\cap L(v)}]
    =d\cdot \sum_{uv\in F'} \sum_{c\in [\chi]}\mathbb{P}[c\in L(u)\cap L(v)]
    =d\cdot e(F')\cdot \chi\cdot\frac{d^2}{\chi^2}
    =e(F')\cdot\frac{d^3}{\chi}
    \end{align*}
    As $\Lambda q^2\leq\chi\leq 2\Lambda q^2$, and $\frac{mq^2}{2}\leq e(F')\leq 2mq^2$ by~\cref{eq:Bound_DeltaF}, we obtain
    \[\frac{md^3}{4\Lambda}\leq\mathbb{E}|\mc{Z}|\leq \frac{2md^3}{\Lambda}.\]
    For $uv\in E(F')$, let $Y_{uv}:=|L(u)\cap L(v)|$, hence such that $|\mc{Z}|=d\cdot\sum_{uv\in F'}Y_{uv}$. As $L(u),L_(v)$ are chosen independently and uniformly among $d$-subsets of $[\chi]$, $Y_{uv}$ is a hypergeometric random variable,
    \[Y_{uv}\sim{\rm Hypergeometric}(\chi,d,d).\]
    In particular, we obtain that 
    \[{\rm Var}(Y_{uv})=d\cdot \frac{d}{\chi}\cdot \frac{\chi-d}{\chi}\cdot\frac{\chi-d}{\chi-1}\leq \frac{d^2}{\chi}.\]
    Observe that if $e,f$ are two vertex-disjoint edges, then $Y_e$ and $Y_f$ are independent. Furthermore, for $uv,uw\in F'$, by independence of $L(v)$ and $L(w)$, we have that $Y_{uv}$ and $Y_{uw}$ conditioning on $L(u)$ are independent, hence
    \begin{align*}
    \mathbb{E}[Y_{uv}Y_{uw}]
    =\mathbb{E}[\mathbb{E}[Y_{uv}Y_{uw}|L(u)]]
    =\mathbb{E}\big[\mathbb{E}[Y_{uv}|L(u)]\cdot \mathbb{E}[Y_{uw}|L(u)]\big].
    \end{align*}
    Moreover, for fixed $L(u)$, $\mathbb{E}[Y_{uv}|L(u)]=\frac{d^2}{\chi}$, therefore we obtain
    \[\mathbb{E}[Y_{uv}Y_{uw}]= \parens*{\frac{d^2}{\chi}}^2 = \mathbb{E}[Y_{uv}]\cdot \mathbb{E}[Y_{uw}].\]
    It follows that, for every two edges $e,f\in F'$,
    \[{\rm Cov}(Y_e,Y_f)=0,\]
    and 
    \[{\rm Var}|\mc{Z}|=d^2\cdot \sum_{uv\in F'}{\rm Var}(Y_{uv})\leq d^2\cdot e(F')\cdot \frac{d^2}{\chi}=d\cdot \mathbb{E}|\mc{Z}|.\]
    By Chebyshev inequality, we obtain
    \begin{align*}
        \mathbb{P}\brackets*{~\big| |\mc{Z}|-\mathbb{E}|\mc{Z}| \big| \geq \frac12\cdot\mathbb{E}|\mc{Z}|~}    
        \leq \frac{4\cdot {\rm Var}|\mc{Z}|}{(\mathbb{E}|\mc{Z}|)^2} = \frac{4d}{\mathbb{E}|\mc{Z}|}
        \leq\frac{16\Lambda}{md^2}\leq \frac{1}{10},
    \end{align*}
    for $q$ large enough. We obtain that, with probability at least $9/10$ we have
    \[\frac{md^3}{8\Lambda}\leq|\mc{Z}|\leq \frac{4md^3}{\Lambda},\]
    as desired.
\end{proof}

\end{document}